\documentstyle[12pt,amssymb]{article}

\input epsf.tex
\begin{document}

\title{Deformations of algebras over operads and Deligne's conjecture}

\author {Maxim Kontsevich and Yan Soibelman}
\maketitle

\newtheorem{thm}{Theorem}
\newtheorem{lmm}{Lemma}
\newtheorem{dfn}{Definition}
\newtheorem{rmk}{Remark}
\newtheorem{prp}{Proposition}
\newtheorem{conj}{Conjecture}
\newtheorem{exa}{Example}
\newtheorem{cor}{Corollary}
\newtheorem{que}{Question}

\vspace{15mm}

{\bf Introduction}

\vspace{5mm}

The deformation theory of associative algebras is a guide for 
developing the deformation
theory of many algebraic  structures. 
Conversely, all the concepts
of what should be the ``deformation theory of everything" must
 be tested in the case
of associative algebras. 

An associative algebra is an algebra over an operad. 
This fact, along with the observation
that in many  examples we are dealing 
with algebras over operads, ``explains" the universality
of associative algebras. It also suggests how to 
develop the deformation theory of algebras
over operads.  
This theory is one of the topics of the present paper. 

Another remarkable fact is the  relationship of the deformation theory
of associative algebras to the geometry of 
configuration spaces of points on surfaces.
One of its incarnations is  Deligne's conjecture. 
A brief history of the conjecture as well as 
its generalizations can be found in [Ko3].
Deligne's conjecture is the second topic of our paper.

The theme, which  motivated the whole project, is 
the Grothendieck-Teichm\"uller 
group ($GT$ for short) and its role in the deformation theory. 
The Grothendieck-Teichm\"uller group can be defined as the automorphism 
group of the tower of the pro-nilpotent completions  of the pure braid groups
(see [Dr]).
The pure braid group of $n$ strings is  the fundamental
group of the configuration space of $n$ points in the plane.
 Deligne's conjecture
shows its relation to the Hochschild complex. 
More evidence for that relation has been found
in deformation quantization. The coefficients in the explicit formulas for
the deformed multiplication given in [Ko1] are periods of algebraic varieties
defined over the field of rational numbers. 
Through the hypothetical
relation to mixed Tate motives 
this fact leads to some conjectures about the action
of the Grothendieck-Teichm\"uller group on the moduli space of
 quantized algebras of functions on a manifold (see [Ko3]).
Deligne's conjecture says that the Hochschild complex
$C^{\cdot}(A,A)$ of an associative (or more general,
$A_{\infty}$) algebra $A$ carries a structure of $2$-algebra
(i.e. an algebra over the operad of chains of the little dics
operad). It follows from our results that the Grothendieck-Teichm\"uller group 
acts (homotopically)
on the moduli space of structures of $2$-algebras on $C^{\cdot}(A,A)$. 
This action is closely related to the action of the 
motivic Galois group described in [Ko3]. We hope to discuss this
topic in detail elsewhere.

The paper is organized in the following way.

Section 1 is devoted to the review of operads.
 It includes a brief introduction
to operads, language of trees, polynomial functors etc.

The deformation theory of algebras over free operads is discussed
in Section 2.
 The ``general" deformation
theory is based on the notion of a formal pointed dg-manifold 
(a formal ${\bf Z}$-graded manifold
with a vector field $d$ of degree $+1$ such that $[d,d]=0$).
We do not use any deep results concerning dg-manifolds.
All necessary facts and definitions can be found in [Ko1]
(note that in [Ko1], [Ko2] formal pointed
dg-manifolds were called formal pointed $Q$-manifolds).
We explain how
to construct a formal pointed dg-manifold which controlls
the deformation theory of an algebra over an operad.
In traditional language this corresponds to a construction
of the deformation functor as a functor of the category
of local Artinian
rings. We explain how this approach leads to the
(homotopy) action of the Grothendieck-Teichm\"uller group
on the dg-manifold conrolling the deformation theory
of the Hochschild complex of an $A_{\infty}$-algebra.

In  Section 3 we outline the strategy and state the theorem
(proved in Section 6) which explain why the $GT$ group 
appears in the deformation
theory of associative (and more generally $A_{\infty}$) algebras.

Our approach is based on the notion of a free resolution of
an operad. This is the subject of Section 4.
 We construct canonical free resolutions via an approach
similar to the one developed by Boardman and Vogt in [BV].
We use free resolutions of operads instead of the  conventional
approach which uses free resolutions of algebras.
This approach to the deformation theory of algebras over operads
is not used very often (although see [M]).

In Section 5 we construct an operad $M$ 
 which acts on the Hochschild complex
of an $A_{\infty}$-algebra. It turns out that $M$ is closely
related to the compactifications of configuration spaces of point
introduced by Fulton and Macpherson in [FM].

Section 6 is devoted to the proof of the theorem from  Section 3.

Section 7 is devoted to  Deligne's conjecture
and its proof.  The proof presented
in this paper is based on the general deformation theory
developed in the previous sections. The strategy is explained in detail
in Section 7. We also suggest certain generalizations
of the original Deligne's conjecture as well as some 
conjectures about the cell structure of the spaces which appear
in the course of its study.
It seems the proof admits a generalization to the higher-dimensional
version of  Deligne's conjecture proposed in [Ko3].

Other approaches to the original
Deligne's conjecture were proposed in [MS], [T], [V].

The Appendix is devoted to a theory of piecewise algebraic chains.
It is suitable
for real semialgebraic manifolds with corners. 
A typical example is the compactification of the configuration space of points
in the plane (Fulton-Macpherson compactification).
This theory is useful for the proof of formality of the operad
of chains of the little discs operad (in the approach of [Ko3]). It is also
useful (but not necessary)
in the proof of  Deligne's conjecture.
 The
usual techniques of simplicial homology will do the job. At the same time we 
feel that the theory of piecewise algebraic chains is appropriate
to the nature of the topic. For this reason we have decided
to include it in the paper. 

{\it Acknowledgements} Second author acknowledges the
support from the Clay Mathematical Institute and IHES.
He thanks IHES for hospitality and stimulating research atmosphere.

\section{Generalities on operads}

\subsection{Polynomial functors, operads, algebras}

The material of this subsection is very well-known (see for ex. [GiKa],
[GeJ], [GeKa], [Ma]). We present it here for completeness and in order to fix
the notation.

Let $k$ be a field of characteristic zero. All vector spaces below
will be $k$-vector spaces unless we say otherwise.
 
We fix a category ${\cal C}$ which is assumed to be
 $k$-linear abelian symmetric monoidal and
 closed under infinite sums and products. We will also assume that it has
 internal $Hom's$. Our main examples will be the category of
 $k$-vector spaces, the category $Vect_{\bf Z}$ of ${\bf Z}$-graded
  vector spaces (with Koszul rule of signs), and the category
  of complexes of $k$-vector spaces.

Suppose we have  a collection of 
representations $F=(F_n)_{n \ge 0}$ of the symmetric groups 
$S_n, n=0,1,...$ in ${\cal C}$
(i.e. we have a sequence of objects $F_n$ together with an action of the group 
$S_n$ on $F_n$ for each $n$).

\begin{dfn}
A polynomial functor $F: {\cal C} \rightarrow {\cal C}$ is defined 
on objects by the formula

$$
F(V)=\oplus_{n \ge 0}(F_n \otimes V^{\otimes n})_{S_n}
$$
where for a group $H$ and an $H$-module $W$ we denote by
$W_H$ the space of coinvariants. Functor
$F$ is defined on morphisms in an obvious way.

\end{dfn}
Notice that having a sequence $F_n$ as above we can define $F_I$ for
any finite set $I$ using isomorphisms of $I$ with the standard set
$\{1,...,|I|\}$, where $|I|$ is the cardinality of $I$. Thus 
$F_{\{1,...,n\}}=F_n$. Technically speaking, we consider a functor $\Phi$
from the groupoid of finite sets (morphisms are bijections)
to the symmetric monoidal category ${\cal C}$. Then we set 
$F_I=\Phi(I)$.

Polynomial functors on ${\cal C}$ form a category ${\cal PF}$ 
if we define morphisms
between two such functors $F$ and $G$ as a vector space of $S_n$-intertwiners
$$
Hom(F,G)=\prod_{n=0}^{\infty}Hom_{S_n}(F_n,G_n)
$$
There is a composition operation $\circ$ on polynomial functors such that
$(F\circ G)(V)$ is naturally isomorphic to $F(G(V))$ for any $V \in {\cal C}$.
We also have a polynomial functor ${\bf 1}$ such that ${\bf 1}_1=1_{\cal C}$
and ${\bf 1}_n=0$ for all $n\ne 1$. Here $1_{\cal C}$ is the unit object in
the monoidal category ${\cal C}$. It is easy to see that in this way we
get a monoidal structure on ${\cal PF}$.

\begin{dfn}
An operad in ${\cal C}$  is a 
monoid in the monoidal category ${\cal PF}$.
In other words it is a polynomial functor $R \in {\cal PF}$ together
with morphisms $m: R\circ R \rightarrow R$ and $u: {\bf 1}\rightarrow R$
satisfying the associativity and the unit axioms.

\end{dfn}

To shorten the notation we will denote the operad $(R,m,u)$ simply
by $R$.
An operad $R$ gives rise to a triple in the category ${\cal C}$.
There is the notion of an algebra over a triple in a category. Hence we can 
use it in order to give a definition of an algebra in ${\cal C}$
over the operad $R$. It is given by an object $V\in {\cal C}$ and a
morphism $R(V) \rightarrow V$ satisfying natural properties of
compatibility with the structure of a triple. 
Equivalently, $V$ is an $R$-algebra iff
 there is a morphism of operads $R \rightarrow  {\cal E}nd(V)$,
where ${\cal E}nd(V)$ is the endomorphism operad of $V$ defined by
$({\cal E}nd(V))_n=\underline{Hom}(V^{\otimes n}, V)$, and $\underline{Hom}$
denotes the internal $Hom$ in ${\cal C}$. The category of $R$-algebras
will be denoted by $R-alg$. There are two
adjoint functors $Forget_R: R-alg\to {\cal C}$ 
and $Free_R:{\cal C}\to R-alg$ such that $Forget_R \circ Free_R=R$.

\begin{dfn} For $X\in Ob({\cal C})$ we call $Free_R(X)$ the free
$R$-algebra generated by $X$.

\end{dfn}

We remind to the reader that there are operads $As$, $Lie$, $Comm$
such that the algebras over them  in the category of vector spaces
are associative, Lie and commutative algebras correspondingly.

\subsection{Colored operads}

There is a generalization of the notion of operad. It is useful 
in order to describe in operadic terms pairs 
{\it (associative algebra A, A-module)},
homomorphisms of algebras over operads, etc.

Let $I$ be set. We consider the category ${\cal C}^I$ consisting of 
families $(V_i)_{i\in I}$ of objects of ${\cal C}$.
 
A polynomial functor $F: {\cal C}^I\rightarrow {\cal C}^I$ is defined by the
following formula:
$$
(F((V_i)_{i\in I}))_j=\oplus_{a:I\rightarrow {\bf Z_{\ge 0}}}
F_{a,j}\otimes_{\prod_iS_{a(i)}}\otimes_{i\in I} (V_i^{\otimes a(i)})
$$
where $a: I\to {\bf Z}_+$ is a map with the finite support, and
$F_{a,j}$ is a representation in ${\cal C}$ of the 
group $\prod_{i\in I}S_{a(i)}$.

Polynomial functors in ${\cal C}^I$ 
form a monoidal category with the tensor product
given by the composition of functors.

\begin{dfn}
A colored operad
is a monoid in this category. 
 
 \end{dfn}
 
Similarly to the case of usual operads
it defines a triple in the category ${\cal C}^I$. Therefore we have
the notion of an algebra over a colored operad.

There exists a colored operad ${\cal OP}$ such that the category of 
${\cal OP}$-algebras is equivalent to the category of operads.

Namely, let us consider the forgetful functor $Operads\to {\cal PF}$.
It has a left adjoint functor. Thus we have a triple in ${\cal PF}$.
As we have noticed before, the category ${\cal PF}$ can be described as
a category of sequences $(P_n)_{n\ge 0}$ of $S_n$-modules.
Then using the representation theory of symmetric groups, we conclude
that the category ${\cal PF}$ is equivalent to the category ${\cal C}^{I_0}$,
where $I_0$ is the set of all Young diagrams (partitions).
Hence a polynomial functor $F:{\cal PF}\to {\cal PF}$ can be described
as a collection $F_{(m_i),n}$ of the representations of the groups
$S_{n,(m_k)}:=S_n\times 
\prod_{k\ge 0}(S_{m_k}\triangleright S_k^{m_k})$,
where $\triangleright$ denotes the semidirect product of groups.

Having these data we can express  any polynomial functor
$F$ on ${\cal PF}$ by the formula:

$$
(F((U_k)_{k\ge 0}))_n=\bigoplus_{(m_k)}F_{(m_k),n}
\otimes_{S_{1,(m_k)}}\bigotimes_{k\ge 0}( U_k^{\otimes m_k})
$$

In particular, one has a functor ${\cal OP}:{\cal PF}\to {\cal PF}$,
which is the composition of the forgetful functor $Operads\to {\cal PF}$
with its adjoint. 
It gives rise to an $I_0$-colored operad ${\cal {OP}}=
({\cal {OP}}_{(m_i),n})$. We will describe it explicitly
in the subsection devoted to trees.

\subsection{Non-linear operads}

We remark that operads and algebras over operads
can be defined for any symmetric monoidal
category ${\cal C}$, not necessarily $k$-linear. In particular, we are 
going to use operads in the categories of sets, topological spaces, etc.

Namely, an operad in ${\cal C}$ is a collection $(F_n)_{n\ge 0}$
of objects in ${\cal C}$, each equipped with an $S_n$-action,
as well as composition maps:

$$F_n\otimes F_{k_1}\otimes...\otimes F_{k_n}\to F_{k_1+...+k_n}$$
for any $n\ge 0,k_1,...,k_n\ge 0$.
Another datum is the unit, which is a morphism ${\bf 1}_{\cal C}\to F_1$.
All the data are required to satisfy certain axioms (see [Ma]).
Analogously one describes colored operads and algebras over operads.
Notice that in this framework one cannot speak about polynomial
functors and free algebras.

This approach has some advantages and drawbacks
(like description of analytic functions  in terms of Taylor
series vs their description in terms of Taylor coefficients).

\subsection {Pseudo-tensor categories}

The notion of colored operad  has been rediscovered 
many times.  
In [BD] the notion of {\it pseudo-tensor category}
was introduced as a generalization of the notion
of symmetric monoidal (=tensor) category. 
The terminology stresses the similarity of operads with tensor
categories. Similar notion was introduced in [B]
under the name {\it multi-linear category}.

This notion is essentially equivalent to the
notion of colored operad. We recall it below,
using the name suggested in [BD].

\begin{dfn} A pseudo-tensor category is given by the following data:

1. A class ${\cal A}$ called the class of objects, and a symmetric
monoidal category ${\cal V}$ called the category of operations.

2. For every finite set $I$, a family $(X_i)_{i\in I}$ of objects,
and an object $Y$, an object $P_I((X_i),Y)\in {\cal V}$ called the
space of operations from $(X_i)_{i\in I}$ to $Y$.

3. For any map of finite sets $\pi: J\to I$ , two families of
objects $(Y_i)_{i\in I}, (X_j)_{j\in J}$ and an object $Z$, a morphism
in ${\cal V}$

$$ P_I((Y_i),Z)\otimes (\otimes_iP_{\pi^{-1}(i)}((X_{j_i}),Y_i))\to
P_J((X_j),Z)
$$
called composition of operations. Here we denote by $\otimes$
the tensor product in ${\cal M}$.

4. For an $1$-element set $\cdot$ and an object $X$, a unit morphism
${\bf 1}_{\cal V}\to P_{\cdot}((X),X)$.

These data are required to satisfy natural conditions. In particular,
compositions of operations are associative with respect to morphisms
of finite sets, and the unit morphisms satisfy the properties analogous
to those of  the identity morphisms (see [BD] or [So] for details).

\end{dfn}

If ${\cal A}$ is a set, then a pseudo-tensor category is exactly the same
as an ${\cal A}$-colored operad in the tensor category ${\cal V}$.

If we take ${\cal V}$ to be the category of sets, and take $I$ above
to be $1$-element sets only, we obtain a 
category with the class of objects equal to ${\cal A}$.

Colored operad  with one color gives rise to an ordinary operad.
 A symmetric monoidal category ${\cal A}$ produces the colored 
operad with $P_I((X_i),Y)=Hom_{{\cal A}}(\otimes_iX_i,Y)$. 

The notion of pseudo-tensor category admits a generalization
to the  case when no action of symmetric group is assumed.
This means that we consider {\it sequences} of objects instead
of {\it families} (see [So]). The new notion generalizes monoidal
categories. In terms of the next subsection this would mean
that one uses planar trees instead of all trees.
One can make one step further generalizing
braided categories.  This leads to colored braided operads
(or pseudo-braided categories). In this case trees in ${\bf R}^3$
should be used.
The deformation theory of such structures can be developed along
the lines of present paper. It will be explained in detail
elsewhere.

\subsection{Trees}
 
\begin{dfn}
A tree $T$  is defined by the following data:

1) a finite set $V(T)$ whose elements are called vertices;

2) a distinguished element $root_T \in V(T)$ called root vertex;

3) subsets $V_i(T)$ and $V_t(T)$ of $V(T)\setminus \{root_T\}$ called the
set of internal vertices and the set of tails
respectively. Their elements are called internal and tail vertices
respectively;

4) a map $N=N_T: V(T) \rightarrow V(T)$.

These data are required to satisfy the following properties:

a) $V(T)= \{root_T\} \sqcup V_i(T) \sqcup V_t(T)$;

b) $N_T(root_T)=root_T$, and $N_T^k(v)=root_T$ for all $v\in V(T)$
and $k \gg 1$;

c) $N_T(V(T))\cap V_t(T)=\emptyset$ ;

d) there exists a unique vertex $v\in V(T), v\ne root_T$ such that
$N_T(v)=root_T$.

\end{dfn}

We denote by $|v|$ the $valency$ of a vertex $v$, which we understand as
the cardinality of the set $N_T^{-1}(v)$.

We call the pairs $(v, N(v))$ $edges$ in the case if $v\ne root_T$. If both
elements of the pair belong to $V_i(T)$ we call the corresponding edge
$internal$. The only edge $e_r$ defined
by the condition d) above is  called the $root$ $edge$. All  edges
of the type $(v,N(v)),v\in V_t(T)$ are called $tail$ $edges$.
We use the notation $E_i(T)$ and $E_t(T)$ for the sets of
 internal and tail edges
respectively. We have a decomposition of the
 set of all edges $E(T)= E_i(T)\sqcup (E_t(T)\cup \{e_r\})$.
 There is a unique tree $T_e$ such that $|V_t(T_e)|=1$ and $|V_i(T_e)|=0$.
  It has  the only tail edge which is also the root edge.

A $numbered$ tree with $n$ tails is by definition a tree $T$ together with
a bijection of sets $\{1,...,n\} \rightarrow V_t(T)$.
We can picture trees as follows

\vskip 1cm
\centerline {\epsfbox {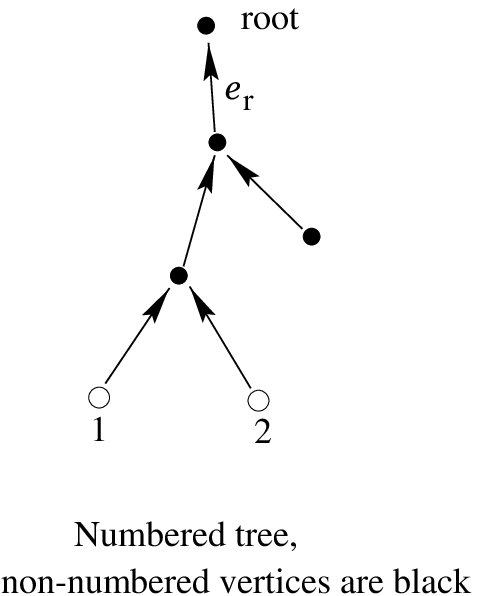} }
\vskip 1cm

Let $R$ be an operad. Any tree $T$ gives a natural way to
compose elements of $R$, 
$comp_T:\otimes_{i\in V_i(T)} R_{N^{-1}(v)}\to R_{V_t(T)}$.

Let us return to the colored operad ${\cal OP}$ and give 
its description using the language of trees.

Namely, ${\cal {OP}}_{(m_i),n}$ is a $k$-vector space generated by
the isomorphism classes of trees $T$ such that:

a) $T$ has $n$ tails numbered from $1$ to $n$;

b) $T$ has $\sum_i m_i$ internal vertices all numbered in such a way
that first $m_0$ vertices have valency $0$, and they are numbered from
$1$ to $m_0$, next $m_1$ internal vertices have valency $1$, and they
are numbered from $1$ to $m_1$, and so on;

c) for every internal vertex $v\in V_i(T)$ the set of incoming edges
$N_T^{-1}(v)$ is also numbered.

An action of the group $S_{(m_k),n}$ is defined naturally:
the factor $S_n$ permutes numbered tails, the factor $S_{m_k}$
permutes numbered internal vertices and the factor $S_k^{m_k}$
permutes  their incoming edges numbered from $1$ to $k$.

The composition is given by the procedure of inserting of a tree 
into an internal vertex of another one. The new numeration is clear.
 We leave these details as well as the proof of the following
 proposition to the reader.
 
\begin{prp} The category of ${\cal OP}$-algebras is equivalent to
the category of $k$-linear operads.

\end{prp}

Let $F$ be a polynomial functor on ${\cal C}$ (see Section 1.1).
 Let us consider a category 
${\cal C}_F$ objects of which are pairs $(V, \phi: F(V) \rightarrow V)$
where $V$ is an object of ${\cal C}$ and $\phi$ is a morphism in ${\cal C}$.
Morphisms of pairs are defined in the natural way. The following lemma
is easy to prove.

\begin{lmm}
The category ${\cal C}_F$ is equivalent to the category of 
$Free_{\cal OP}(F)$-algebras. 

\end{lmm}

We call $P=Free_{\cal OP}(F)$ the $free$ $operad$ generated by $F$.

Components $P_n$ of the functor $P$ can be defined explicitly as follows.
 
 Let $Tree(n)$ denotes
the groupoid of numbered trees with $n$ tails, $|Tree(n)|$ denotes the 
set of classes of
isomorphisms of these trees . We denote the class of isomorphism
of $T$ by $\lbrack T \rbrack$.
 Then we have

$$P_n=Free_{{\cal OP}}(F)_n=\oplus_{ \lbrack T \rbrack\in |Tree(n)|}
(\otimes_{v\in V_i(T)}F_{N^{-1}(v)})_{AutT}
$$

\section{Deformations and differentials in  free operads}

Let $F$ be a polynomial functor, $P=Free_{\cal OP}(F)$
 be the corresponding free operad.
Let $g_P$ be the Lie algebra (in the symmetric monoidal category ${\cal C}$)
of derivations of the operad $P$.
Then, as an object of ${\cal C}$:

$$
g_P=\prod_{n \ge 0}\underline{Hom}(F_n,P_n)^{S_n}
$$
where $W^H$ denotes the space of $H$-invariants of an $H$-module $W$
and $\underline{Hom}$ denotes the internal $Hom$ in ${\cal C}$.
This observation follows from the fact that 
$Hom_{\cal PF}(F,Forget_{\cal OP}(G))=Hom_{{\cal OP}-alg}(Free_{\cal OP}(F),G)$.

>From now on we suppose that ${\cal C}$ is the category $Vect_{\bf Z}$
of ${\bf Z}$-graded vector spaces. Then $g_P$ is a graded Lie algebra
with the graded components $g_P^n$.

\begin{dfn}
A structure of a
differential-graded operad on $P$ which is
 free as a graded operad is given by 
an element $d_P \in  g_P^1$ such that $\lbrack d_P,d_P \rbrack=0$.
\end{dfn}

 The definition  means
 that $P$ can be considered as an operad in the symmetric monoidal category
of complexes, and it is free as an operad in the category $Vect_{\bf Z}$.
Sometimes we will denote the 
corresponding operad in the category of complexes by $\widehat P$.

One of our purposes will be to use $\widehat P$ for constructing resolutions
of dg-operads, and subsequently the deformation theory of  algebras 
over them.

\begin{dfn}
A dg-algebra over $(P,d_P)$ (or simply over $P$)  is an algebra
over $\widehat P$ in the category of complexes.
\end{dfn}

Notice that the deformation theory of the pair $(P,d_P)$ is the same as
the deformation theory of $d_P$ (since $P$ is free and therefore rigid). 

\begin{dfn}
The formal pointed dg-manifold
associated with the
differential graded Lie algbera
 $(g_P, \lbrack d_P, \cdot \rbrack)$ controls the deformation
theory of $(P,d_P)$.
\end{dfn}

Now we are going to describe the deformation theory of dg-algebras
over $P$ ($\widehat P$-algebras). 
In what follows we will
often speak about
points of ${\bf Z}$-graded manifolds. It will always mean
$\Lambda$-points, where $\Lambda$  is an auxiliary ${\bf Z}$-graded
commutative associative algebra (in general without the unit).
In the case of formal manifolds we take $\Lambda$ to be nilpotent.

Let us describe the formal pointed dg-manifold controlling deformations
of a $P$-algebra $V$. 
We have the following  graded vector space

$$
{\cal M}=
{\cal M}(P,V)=(\underline{Hom}(V,V))\lbrack1\rbrack \oplus
\underline{Hom}(F(V),V)
$$

We denote by ${\cal M}^n, n\in {\bf Z}$ the graded components of ${\cal M}$.

The structures of a complex on a graded vector space $V$ and an action
of $P$ on $V$ define a point $(d_V, \rho) \in {\cal M}^0=
Hom_{Vect_{\bf Z}}(k,{\cal M})$. We consider here $d_V$ and $\rho$
as morphisms of graded vector spaces.
The equation $d_V^2=0$ and the condition
of compatibility of $d_V$ and $\rho$ can be written in the form
$d_{\cal M}(d_V,\rho)=0$, where 
$d_{\cal M}(d_V,\rho)=(d_V^2,\xi(d_V,\rho))\in {\cal M}^1$,
for some $\xi(d_V,\rho)\in \underline{Hom}(F(V),V)$.
It is easy to see that the assignment $(d_V,\rho)\mapsto d_{\cal M}(d_V,\rho)$
defines an odd vector field $d_{\cal M}$ on the 
``infinite-dimensional graded manifold" ${\cal M}$. 
A zero of this vector field corresponds to
a structure of a complex  on $V$ together
with a compatible structure of a dg-algebra over $ P$.
This gives a bijection between the set of zeros and the set
of such structures.

It is easy to check that $\lbrack d_{\cal M},d_{\cal M} \rbrack=0$.
Therefore a formal neighborhood of a fixed point $(d_V,\rho)$
of $d_{\cal M}$
becomes a formal pointed dg-manifold.

\begin{dfn}
The deformation theory of a dg-algebra $V$ is
controlled by this formal pointed dg-manifold.
\end{dfn}

\begin{rmk} Operads are algebras over the colored
operad ${\cal OP}$. One can show that the deformation theories
for an ${\cal OP}$-algebra $\widehat P$ described in the last two
definitions are in fact equivalent. 
\end{rmk}

Let $P^{(0)}$ be an operad in $Vect_k$. In order to define
the deformation theory of $P^{(0)}$-algebras, one needs
to choose a resolution $\phi:P\to P^{(0)}$, where
$P=Free_{\cal OP}(F)$ as a ${\bf Z}$-graded operad, and 
$\phi$ is a quasi-isomorphism.
One wants to be sure that the deformation theory does not
depend on the choice of the resolution.
This can be achieved by assuming that (see for ex. [M2])
:

a) $F$ admits a filtration (as a polynomial functor)
$F=\cup_{j\ge 1}F^{(j)}, F^{(j)}\subset F^{(j+1)}$
such that $d_P(F^{(0)})=0$ and 
$d_P( F^{(j)})\subset Free_{\cal OP}(F^{(j-1)}), j\ge 1$;

b) $\phi:P\to P^{(0)}$ is a an epimorphism.

All resolutions used in the paper will satisfy these properties.

\subsection{Example: $A_{\infty}$-operad and $A_{\infty}$-algebras}

Let $V\in Vect_{\bf Z}$ and $m_n:V^{\otimes n}\to V[n-2], n\ge 2$ 
be a sequence of morphisms. It gives rise to an action on $V$ of the free
operad $P=Free_{\cal OP}(F)$ where

$$F(V)=\oplus_{n\ge 2}V^{\otimes n}[n-2].$$

Then $F_n=k[S_n]m_n\otimes k[1]^{\otimes(n-2)}$. 
This notation means that we consider
$F_n$ as a  space (with the grading shifted by $n-2$)
of the regular representation of the group algebra
of the symmetric group $S_n$. This space is generated by 
an element which we denote by $m_n$.

The differential $d_P\in g_P$ (equivalently, a structure
of a dg-operad on $P$) is defined by the standard formulas:

$$ d_P(m_2)=0,$$

$$ d_P(m_n)(v_1\otimes...\otimes v_n)=
\sum_{k+l=n}\pm m_k(v_1\otimes...\otimes v_i
\otimes m_l(v_{i+1}\otimes...\otimes v_{i+l})\otimes...\otimes v_n), n>2.$$

We do not specify signs in these well-known formulas.
In Section 5 we propose a general framework allowing to fix
signs in the formulas like the one above.

\begin{dfn} The dg-operad ${\cal A}_{\infty}=(P,d_P)$ 
is called the $A_{\infty}$-operad.
Algebras over this dg-operad are called $A_{\infty}$-algebras.

\end{dfn}

Deformations of an $A_{\infty}$-algebra $A$ are controlled by the truncated
Hochschild complex

$$C_+^{\cdot}(A,A)=\prod_{n\ge 1}Hom_{Vect_{\bf Z}}(A^{\otimes n},A)[-n]$$

More precisely,
let $A$ be a graded vector space. We define a graded vector space
of all Hochschild cochains of $A$ as

$$
C^{\cdot}(A,A)=\prod_{n\ge 0}Hom_{Vect_{\bf Z}}(A^{\otimes n},A)
\lbrack -n \rbrack
$$
Then $C^{\cdot}(A,A)[1]$ can be 
equipped with the structure of a graded Lie algebra
 with the Gerstenhaber bracket (the latter appears naturally
 if we interpret  Hochschild cochains as coderivations
 of the free coalgebra $\oplus_{n\ge 0}(A[1])^{\otimes n}$).

Let us consider an element $m=(m_1,m_2...)\in C_+^{\cdot}(A,A)[1]$ 
of degree $+1$
such that $\lbrack m,m\rbrack=0$. Such an element defines a differenitial
$d=m_1$ on $A$, and the sequence $(m_2,m_3,...)$ gives rise to a structure
of an $A_{\infty}$-algebra on $(A,m_1)$.

Then we
can make $C^{\cdot}(A,A)$ into a complex (Hochschild complex) with the
differential $d_m=\lbrack m,\cdot\rbrack$. It is easy to see
that in this way we get a differential graded Lie algebra
(DGLA for short) $(C^{\cdot}(A,A)[1],d_m)$. The truncated
Hochschild complex $C_+^{\cdot}(A,A)[1]$ is a DGLA subalgebra.
According to the general theory (see [Ko1]) both DGLAs define formal
pointed dg-manifolds, and therefore give rise to the
deformation functors. This is a straightforward generalization
of the well-known deformation theory of associative algebras.

In a sense the full Hochschild complex controls deformations of the
 $A_{\infty}$-category with one object, such that
 its endomorphism space is equal to $A$. 

The deformation theory of $A_{\infty}$-categories is not in the scope
of present paper. It will be explained elsewhere.
Nevertheless we will refer to the formal dg-manifold associated with
$C^{\cdot}(A,A)[1]$ as to the  {\it moduli space of $A_{\infty}$-categories}.
Similarly, the formal dg-manifold associated with $C_+^{\cdot}(A,A)[1]$
will be called the {\it moduli space of $A_{\infty}$-algebras}.
(All the terminology assumes that we deform a given $A_{\infty}$-algebra $A$).

The moduli space of $A_{\infty}$-algebras is the same as
 ${\cal M}({\cal A}_{\infty},A)$ in the previous notation.
 Similarly we will denote the moduli space of $A_{\infty}$-categories
 by ${\cal M}_{cat}({\cal A}_{\infty},A)$. The natural inclusion
 of DGLAs $C_+^{\cdot}(A,A)[1]\to C^{\cdot}(A,A)[1]$ induces
 a dg-map 
 ${\cal M}({\cal A}_{\infty},A)\to {\cal M}_{cat}({\cal A}_{\infty},A)$
 (dg-map is a morphism of dg-manifolds).

Let us remark that the operad ${\cal A}_{\infty}$ gives rise
to a free resolution of the operad ${ As}$. Algebras
over the latter are associative algebras without the unit.

\begin{rmk}
It is interesting to describe deformation theories of
free resolutions of the classical operads ${ As}$,
${Lie}$, ${ Comm}$. It seems that for an arbitrary 
free resolution $P$ of either of these operads the following
is true: $H^i(g_P)=0$ for $i\ne 0$, $H^0(g_P)=k$.
This one-dimensional vector space gives rise to 
the rescaling of operations, like $m_n\mapsto \lambda^n m_n$
in the case of $A_{\infty}$-algebras.

\end{rmk}

\subsection{Homotopical actions of Lie algebras}

Let $g$ be a Lie algebra acting on a formal dg-manifold $(Y,d_Y)$.
This means that we have a homomorphism of Lie algebras
$g \rightarrow Der(Y), \gamma \mapsto \hat\gamma$ 
where $Der(Y)$ is the Lie algebra of vector fields
on $Y$ preserving ${\bf Z}$-grading an $d_Y$.

We can make $Z=Y\times g\lbrack 1\rbrack$ 
into a formal dg-manifold introducing
an odd vector field by the following formula

$$
d_Z(y,\gamma)=(d_Y(y)+\hat \gamma, \lbrack \gamma,\gamma \rbrack/2)
$$

Then $\lbrack d_Z,d_Z \rbrack=0$. We can make $g\lbrack 1\rbrack$ into
a formal dg-manifold using the odd vector field $d_{g\lbrack 1\rbrack}$  
arising from the
Lie bracket. Then the natural projection $(Z,d_Z) \rightarrow 
(g\lbrack 1 \rbrack,d_{g\lbrack 1\rbrack})$ becomes a
dg-bundle 
(cf. [Ko2]).

This contsruction of a dg-bundle out of a group or Lie
algebra acting on a dg-manifold motivates the following definition.

\begin{dfn} Let $g$ be a  Lie algebra.
Homotopical $g$-action on a formal dg-manifold $(Y, d_Y)$ is 
a dg-bundle $\pi: (Z,d_Z) \rightarrow (g\lbrack 1 \rbrack,
d_{g\lbrack 1\rbrack})$ together with an isomorphism of dg-manifolds 
$(\pi^{-1}(0),d_Z) \simeq (Y,d_Y)$.

\end{dfn}

\begin{rmk}
It was pointed out in [Ko2] that 
in this case $g$ acts on the cohomology of all
complexes naturally associated with $(Y,d_Y)$ ( like the tangent space
at a zero point of $d_Y$, the space of formal functions on $Y$, etc.).

\end{rmk}

We remark also that if $g$ is a DGLA then the same definition can be given.
It can be also generalized to the case when the base of a equivariant dg-bundle
is an arbitrary dg-manifold with a marked $d$-stable point.

Suppose that $F$ is a polynomial functor in the category of ${\bf Z}$-graded
vector spaces, $P=Free(F), V \in Vect_{\bf Z}$. We apply the general
scheme outlined above to the case $Y={\cal M}(P,V), g=g_P$.
Obviously $g$ acts on the dg-manifold $\underline{Hom}(F(V),V)$,
equipped with the trivial odd vector field.

Let us consider the graded vector space
$$
{\cal N}= \underline{Hom}(V,V)[1]\oplus \underline{Hom}(F(V),V)
\oplus g_P\lbrack 1 \rbrack
$$

Let $d_V\in \underline{Hom}(V,V)[1]$ makes $V$ into a complex,
$\gamma=d_P\in g_P[1]$ satisfies the equation $[d_P,d_P]=0$ and
$\rho\in \underline{Hom}(F(V),V)$ makes $V$ into a dg-algebra
over $(P,d_P)$.

We consider the formal neighborhood of the point $(d_V,\rho,d_P)$
in ${\cal N}$, and define an odd vector field by the formula

$$d_{\cal N}(d_V,\rho,d_P)=(d_V^2,\xi(d_V,\rho)+\hat {d}_P,[d_P,d_P]/2)$$

The notation here is compatible with the one for ${\cal M}$.

One can check that $[d_{\cal N},d_{\cal N}]=0$. Thus the formal neighborhood
becomes a formal dg-manifold. It controls deformations
of pairs {\it (an operad, an algebra over this operad)}.

The natural
projection $\pi: {\cal N} \rightarrow g_P\lbrack 1 \rbrack$ is a 
morphism of formal dg-manifolds . Here on $g_P\lbrack 1 \rbrack$ we use
the odd vector field $d_{g_P \lbrack 1\rbrack}$ defined by the Lie
bracket. Then the formal scheme of zeros
 of $d_{g_P\lbrack 1 \rbrack}$ corresponds
to the structures of a dg-operad on $P$. The fiber over a fixed point
$x\in g_P[1]$ is a dg-manifold with the differential induced from ${\cal N}$.
Then the formal neighborhood of a fixed point in $\pi^{-1}(x)$ controls
deformations of $\widehat P$-algebras. 

We conclude that the Lie algebra 
of derivations of an operad acts homotopically
on the moduli space of algebras over this operad. 

\section{Hochschild complex and operads}

 This section serves as a sort of 
a ``second introduction", outlining objectives and 
the strategy of the rest of the paper.

One of our aims will be to construct a dg-operad of the type $\widehat P$
(i.e. free as a graded operad)
acting naturally on the Hochschild complex of an arbitrary $A_{\infty}$-
algebra. 

In Section 5 we are going to construct
an operad $M$ which acts naturally on the full Hochschild complex
$C^{\cdot}(A,A)$ of an $A_{\infty}$-algebra $A$
as well as on $C_+^{\cdot}(A,A)$. There is a natural
free resolution $P$ of the operad $M$, so that $C:=C^{\cdot}(A,A)$ becomes
a $\widehat P$-algebra. Then we can say that there is a
dg-map of the moduli space  of  $A_{\infty}$-categories 
to the moduli space ${\cal M}(P,C)$ of structures of 
$\widehat P$-algebras on the graded vector space $C$.

>From the point of view of deformation theory it is not very natural
to make constructions of the type {\it algebraic structure} $\to$
{\it another algebraic structure} (like our construction
$A_{\infty}-algebras\to M-algebras$). It is more natural
to extend them to morphisms between the formal pointed
dg-manifolds controlling the deformation theories of  structures.

We will construct an explicit dg-map
${\cal M}({\cal A}_{\infty},A)\to {\cal M}(P,C)$ as well as 
a dg-map ${\cal M}_{cat}({\cal A}_{\infty},A)\to {\cal M}(P,C)$,
such that the one is obtained from another by the restriction
from the moduli space of algebras to the moduli space of categories.

The operad ${\cal A}_{\infty}$ is augmented, i.e. equipped with a morphism
of dg-operads $\eta: {\cal A}_{\infty} \rightarrow Free(0)$.
 Here $Free(0)$ is the
trivial operad : $Free(0)_{1}=k, Free(0)_{(n\ne1)}=0$. 
Since $A$ (as any graded vector space) is
an  algebra over $Free(0)$, it becomes  also an algebra 
over ${\cal A}_{\infty}$.
Any structure of an ${\cal A}_{\infty}$-algebra on $A$
 can be considered as a deformation
of this trivial structure. Notice also that in the
previous notation the augmentation morphism defines a point in the
dg-manifold ${\cal M}({\cal A}_{\infty},C)$, where $C=C^{\cdot}(A,A)$.
Therefore it is sufficient to work in the formal
neighborhood of this point.

Notice that we can consider also the moduli space of 
 structures of a complex on the graded vector space
$C^{\cdot}(A,A)$ where $A$ is an arbitrary graded vector space.
It gives rise to a formal dg-manifold. 
There are natural morphisms to it from the formal dg-manifold
of the moduli space of $A_{\infty}$-categories and from the
formal dg-manifold of the moduli space of
structures of $\widehat P$-algebras on $C^{\cdot}(A,A)$.
Theorem $1$ below  combines all three morphisms discussed above into
a commutative diagram. Let us make it more precise.

First we formulate a simple general lemma, which will be applied
in the case $V=C[2]$.

Let $V$ be an arbitrary graded vector space, $d_V$ an odd vector field on $V$
(considered as a graded manifold)
such that $\lbrack d_V,d_V\rbrack=0$. Thus we get a dg-manifold.
The graded vector space $H=H(V)=\underline{Hom}(V,V)\lbrack 1\rbrack$ is
a dg-manifold with $d_H(\gamma)=\gamma^2$. To every point $v\in V$
we assign a point in $H$ by taking the first
Taylor coefficient $d^{(1)}_V(v)$ of $d_V$ at $v$.
In this way we obtain a map $\nu: V\rightarrow H$.

\begin {lmm}
The map $\nu$ is a  morphism of dg-manifolds.

\end{lmm}
$Proof$. Let us write in local coordinates $x=(x_1,...,x_n)$ the vector field
 $d_V=\sum_i\phi_i\partial_i$ where $\partial_i$ denotes the partial
 derivative with respect to $x_i$,
 and $\phi_i$ are functions on $V$.
  Then the map $\nu$ assigns to
 a point $x$ the matrix $M=(M_{ij}(x))$ with $M_{ij}=\partial_j\phi_i$.
Then direct computation shows that the condition $\lbrack d_V,d_V \rbrack=0$
implies that the vector field $\dot x=d_V(x)$ is mapped to the vector
field $\dot M=d_H(M)=M^2$. $\blacksquare$

Let $A$ be a graded vector space endowed with the trivial 
$A_{\infty}$-structure, and
$C=C^{\cdot}(A,A)=\prod_{n\ge 0}Hom_{Vect_{\bf Z}}(A^{\otimes n},A)$
 be the graded space of Hochschild
cochains.
Since $C \lbrack 1 \rbrack $ carries a structure of a
graded Lie algebra (with the Gerstenhaber bracket), it gives rise to
the structure of a 
dg-manifold on $C \lbrack 2 \rbrack$, which is the same as
${\cal M}_{cat}({\cal A}_{\infty},A)$. We will denote it by $(X,d_X)$
(or simply by $X$ for short).

Even for the trivial $A_{\infty}$-algebra structure on $A$,
we get a non-trivial $P$-algebra structure on $C$. The corresponding
moduli space ${\cal M}(P,C)$ will be denoted by $(Y,d_Y)$ (or $Y$ for short).

 There is  a natural
morphism of dg-manifolds $p:Y \rightarrow \underline{Hom}
(C,C)\lbrack 1 \rbrack=H$ (projection of $Y={\cal M}(P,C)$ to the first
summand).

The following theorem will be proved later in the paper.

\begin{thm} There exists a $GL(A)$-equivariant
morphism of dg-manifolds $f:X\rightarrow Y$
such that $pf=\nu$.
\end{thm}

Moreover we will  present an explicit construction of
the morphism.

Suppose that $A$ is an $A_{\infty}$-algebra.
Geometrically the structure of an $A_{\infty}$-algebra
on the graded vector space $A$ gives rise to a point $\gamma\in X=C[2],
C=C^{\cdot}(A,A)$ such that $d_X(\gamma)=0$. 
 Indeed, the definition  
can be written as 
$\lbrack \gamma, \gamma \rbrack=0$. Thus we get a differential in $C$
(commutator with $\gamma$) making it into a complex. The structure
of a complex on the graded vector space $C$ gives rise to  
 a zero of the field $d_H$ in the dg-manifold 
$H=\underline {Hom}(C,C)\lbrack 1\rbrack$. Theorem $1$ implies
that $f(\gamma)$ is a zero of the vector field $d_{{\cal M}(P,C)}$.
Therefore the Hochschild complex $(C, \lbrack \gamma,\cdot \rbrack)$ 
carries a structure of a dg-algebra over $P$. 

Our next aim is to uncover the geometric origin of the operad $P$.
It will be related to the configuration space of discs inside
of the unit discs in the plane. More precisely, the
 operad $Chains(E_2)$ of  chains on the 
little  discs operad (see [Ko3] and Section 7 below)
is quasi-isomorphic to $\widehat P$. Here we use either
usual singular chains or piecewise algebraic chains (see Appendix).
In fact we are going to construct explicitly a morphism
$\widehat P\to Chains(E_2)$ which gives the homotopy equivalence
(to be more precise we will do that for the operad $Chains(FM_2)$
which is quasi-isomorphic to $E_2$). Then using the fact that
both dg-operads are free as graded operads, we invert this
quasi-isomorphism. This gives a structure of an
$Chains(E_2)$-algebra on the Hochschild complex of an
$A_{\infty}$-algebra. This result is known as Deligne's
conjecture.

Let us recall that there is a notion of $d$-algebra, $d\in {\bf Z}_+$  
(see for ex. [Ko3]). Namely, a graded vector space $V$ is called
a $d$-algebra if it is an algebra over the operad $Chains(E_d)$
(chains of the topological operad of little $d$-dimensional discs).
 Thus the
moduli space ${\cal M}(P,C)$ can be thought as a moduli space
of structures of a $2$-algebra on a graded vector space $C$.
Then our Theorem $1$ says that there is a $GL(A)$-equivariant
morphism of the moduli space of $A_{\infty}$-categories with
one object to
the moduli space of $2$-algebras.

\begin{rmk} In the unpublished paper [GJ] the name $d$-algebras
was reserved for algebras over the operad $H_{\cdot}(Chains(E_d))$.
It was proved by Tamarkin in [T] and by the first author in [Ko3]
that there exists a (non-canonical) quasi-isomorphism between the
operad $Chains(E_d)$ and its 
homology operad $H_{\cdot}(Chains(E_d))$
(in other words the operad $Chains(E_d)$ is formal).

\end{rmk}

Let $g_P=\underline{Der}P$  means
as before the DGLA of derivations of $\widehat P$. Then $g_P$ acts
on the moduli space of $\widehat P$-algebras. 
On the other hand, from the point of view of deformation
theory, we can replace the operad of little discs $E_2$ by
the operad of configurations of points in the plane
(properly compactified). This is the operad $FM_2$ mentioned
above. There is a natural action of the Grothendieck-Tiechm\"uller
groups on the rational homotopy type of the latter.
It gives rise to a morphism of $L_{\infty}$-algebras
$Lie(GT)\lbrack 1 \rbrack \rightarrow 
(g_P, \lbrack d_P ,\cdot \rbrack)$
where $GT$ is the Grothendieck-Teichm\"uller group.

Therefore one has a homotopical action of the Lie algebra $Lie(GT)$ on
the moduli space of $\widehat P$-algebras.

\section{ Free resolution of a dg-operad }
In this subsection we recall  well-known constructions of 
free resolutions of operads (see e.g. [GJ]).

Let $k$ be a field as before,
 $R$ be a dg-operad  over $k$. The aim of this subsection is to construct
canonically a dg-operad $P_R$ over $k$, which is free as a graded operad, and a
quasi-isomorphism $P_R \rightarrow R$. Then $P_R$ will be a free resolution 
of $R$. In this subsection we will assume that $R$ is non-trivial, which 
means that the unit operation from $R_1$ is not equal to zero.

\subsection{Topological construction}

We will mainly follow [BV].

Let $O= (O_n)_{n\ge 0}$ be a topological operad (i.e.
all $O_n$ are $S_n$-topological spaces and all operadic morphisms are
continuous). We describe (following [BV]) a construction of a topological
operad $B(O)= (B(O)_n)_{n\ge 0}$ together with a morphism of topological
operads $B(O)\rightarrow O$ which is homotopy equivalence.

To simplify the exposition we assume that $S_n$ acts freely on $O_n$ for
all $n$.

Each space $B(O)_n$ will be the quotient of 

$$
\overline{B(O)}_n= \bigsqcup_{\lbrack T\rbrack,  T\in Tree(n)}
(\lbrack 0, +\infty\rbrack^{E_i(T)}\times\prod_{v \in V_i(T)}O_{N^{-1}(v)})/AutT
$$
under the relations described in the following way.

Let us consider the elements of $\overline{B(O)}_n$ as numbered trees with 
elements of $O$ attached to the internal vertices, and the $length$
$l(e)\in \lbrack 0, +\infty\rbrack$ attached to every edge $e$.
We require that all external edges (i.e. root edge and the tail edges)
have lengths $+\infty$.

We impose two type of relations.

1) We can delete every vertex $v$ of valency 1 if it contains the unit
of the operad, replacing it and the attached two edges of lengths $l_i, i=1, 2$
by the edge
with the length $l_1+l_2$ .
 We use here the usual assumption: $l+ \infty= \infty+l=
\infty$.

2) We can contract every internal edge $e=(v_1,v_2), v_2=N(v_1)$ of the length 
$0$
and compose in $O$ the operations attached to $v_i, i=1,2$.

We depict the trees and relations below.

\vskip 1cm
\centerline{\epsfbox {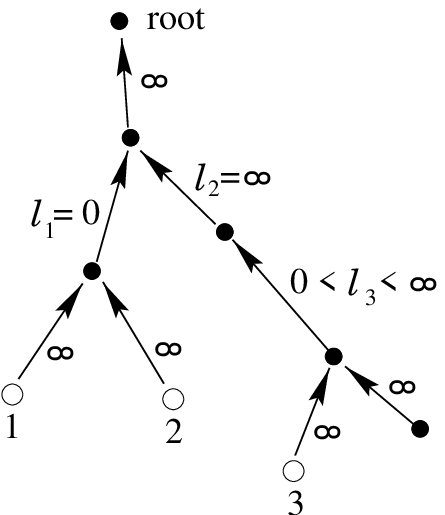}}
\vskip 1cm

\vskip 1cm
\centerline{\epsfbox {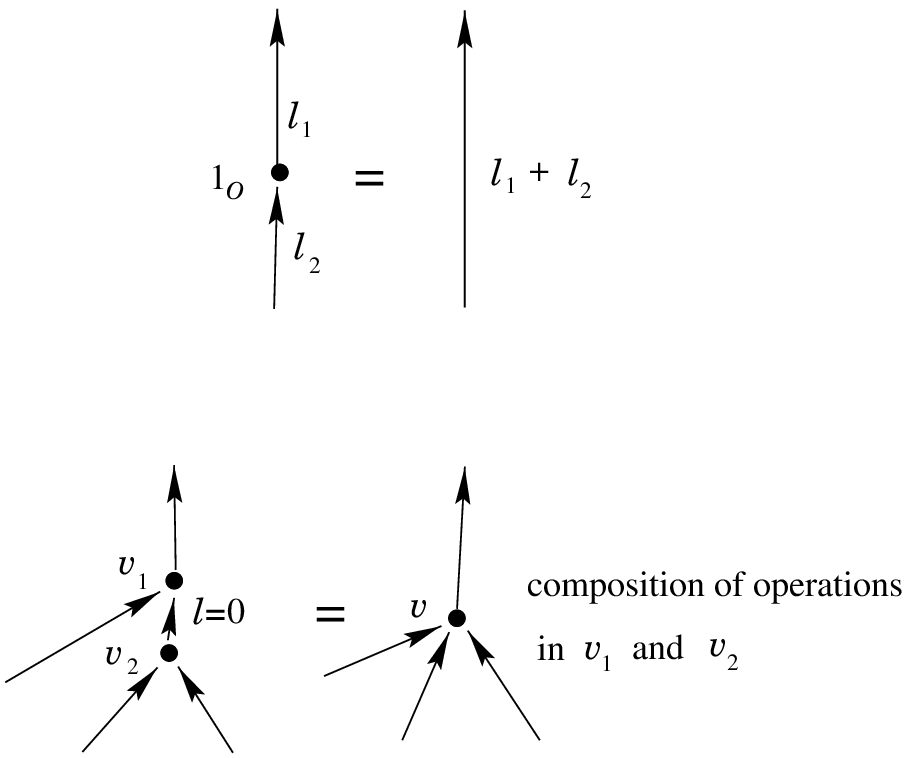}}
\vskip 1cm

Let us describe how $B(O)$ can act naturally
on a topological space.
 
 Let $X$ be a topological space, 
  and $g^t:X\to X, t\in [0,+\infty)$
 a $1$-parametric semigroup of continuous maps acting on $X$.

We assume that the map $[0,+\infty)\times X\to X, (t,x)\mapsto g^tx$
extends continuosly to $[0,+\infty]\times X$.
We denote by $Y$ the image of $g^{\infty}:=lim_{t\to \infty}g^t$.
Clearly the subspace $Y$ is a homotopy retract of $X$.

Suppose that a topological operad $O$ acts on $X$,
 i.e. that we are given continuous maps
$O_n\times X^n \rightarrow X, n \ge 0$, satisfying the usual properties.
We can construct an action of $B(O)$ on $Y$ as follows.
Let $\gamma \in O_n, t, t_i \in {\bf R}_+ \cup \{+\infty\}, 
x_i \in X, 1\le i\le n$. Then we
assign to these data the point $x=g^t\gamma(g^{t_1}x_1,...,g^{t_n}x_n)$
of $X$. We define the composition of such operations in the natural way. 

We can interpret the parameters $t, t_i$ above as lengths of edges of
 trees . Putting $t=+\infty$ we obtain an action
of $B(O)$ on the homotopy retract $Y$.

\subsection{Free resolutions of linear operads}

Let $R$ be a dg-operad over a field $k$.
To  describe its free resolution $P_R$ we need a special
class of trees described below. 

For every $n\ge 0$ we introduce a groupoid ${\cal T}(n)$ of $marked$ $trees$
with $n$ tails. An object of ${\cal T}(n)$ is a numbered tree $T\in Tree(n)$
and a map to a $3$-element set
$l=l_T: E(T)\rightarrow \{0,finite,+\infty\}$ such that
$l_T(\{e_r\}\cup E_t(T))=\{+\infty\}$.
 Notice that in the  case of topological operads
 the component $\overline{B(O)_n}$ is stratified
naturally with the strata labeled by equivalence classes $|{\cal T}(n)|$
of marked trees. The label of an edge $e$ of the corresponding marked tree
is $0$ if $l(e)=0$, is $finite$ if $l(e)\in (0,+\infty)$ and is $+\infty$
if $l(e)=+\infty$. According to this description, we call them {\it zero,
finite, infinite} edges respectively. We denote these sets of edges
by $E_{zero},E_{finite}$ and $E_{infinite}$ correspondingly.

We will give three different but equivalent
descriptions of the  operad $P=P_R$ as a graded
operad. Then we define a differential on $P$.

{\it Description 1. }

Let 
$$
\bar P_n=\bigoplus_{\lbrack T\rbrack ,T \in {\cal T}(n)}
(\otimes_{v\in V(T)}R_{N^{-1}(v)}
\lbrack J_T\rbrack)_{AutT}
$$
where $AutT$ is the group of automorphisms of the tree $T$,
$J_T=l_T^{-1}(finite)$,
and for any graded vector space $W$ and a finite set $J$ we use
the notation $W\lbrack J \rbrack=W\otimes k\lbrack 1\rbrack^{\otimes J}$
(shift of the grading by $J$).

Notice that  the dimension of the corresponding 
stratum of ${\cal T}(n)$ is the cardinality of the set $J_T$, i.e.
the number of finite edges.

Then $(\bar P_n)_{n \ge 0}$ evidently form a graded operad $\bar P$.
It is a $k$-linear analog of the operad $\overline{BO}$.

The operad $\bar P$ contains
a subspace $I_P$ generated by the following relations

1) if the length of an edge $(w,v)$ is $0$ we contract it and make the
composition in $R$ of the operations attached to $w$ and $v$
(cf. description for $B(O)$);

2) if $v$ is a vertex of $T\in \overline P$ such that $|v|=1$
 and $1_R\in R_1$ is attached
to $v$, then $T$ belongs to $I_P$ if the following holds:
lengths of both edges attached to $v$ are non-zero edges, and at
least one of them is finite. 
If both
edges are infinite, we remove the vertex and two attached edges, 
replacing them by an infinite edge.

One can check easily that $I_P$ is a graded ideal in $\bar P$.
We denote by $P$ the quotient operad $\bar P/I_P$.

{\it Description 2.}

We define $P_n$ in the same way, but making the summation
over all trees without edges of zero length. We also drop the relation 1)
from the list of imposed relations (there are no edges with $l=0$ ).

{\it Description 3.}

We define an operad $R^{\prime}$ such as follows:

$R^{\prime}_n=R_n$ for $n \ne 1$, $R^{\prime}_1$ is a
complement to the subspace $k\cdot 1_R$ in $R_1$.

Then we define $P_n$ as in Description 2, but using $R^{\prime}$
instead of $R$ and dropping both relations 1) and 2).

It is clear that this description defines a free graded operad.

Equivalently, it can be described as a free graded operad $P$ such that

$$
P=Free (Cofree^{\prime} (R^{\prime}\lbrack 1 \rbrack))\lbrack -1 \rbrack
$$
Here $Cofree(L)$ means a dg-cooperad generated by $L$ which is cofree
as a graded co-operad, and
${\prime}$ denotes the procedure of taking a
(non-canonical)  complement to the subspace
generated by the unit (or counit in the case of a co-operad)
 as described above in the case of $R$.

In this description the generators of $P$ correspond to such trees $T$ in
${\cal T}=({\cal T}(n))_{n\ge 0}$ that every $T$ has at least one internal 
vertex, all internal edges are finite and there are no zero-edges in $T$.

\begin{prp}
All three descriptions give rise to isomorphic  graded free operads
over $k$.

\end{prp}
$Proof$. Straitforward. $\blacksquare$

We can make $\bar P$ into a dg-operad introducing a differential $d_{\bar P}$.
We use the Description 1 for this purpose.

The differential $d_{\overline P}$ is naturally decomposed into the sum of
 two differentials: 
 
 $d_{\overline P}=\tilde{d}_R+d_{\cal T}$ where

a) the differential $\tilde{d}_R$ arises from the differential $d_R$ in $R$;

b) the differential $d_{\cal T}$ arising from the stratification
of ${\cal T}(n)$: it either contracts a finite edge or makes it 
into an infinite edge .

To be more precise, let us consider the following object $\Delta$
in ${\cal C}= Vect_{\bf Z}$: $\Delta^{-1}=1_{{\cal C}},
\Delta^0=1_{{\cal C}}\oplus 1_{{\cal C}}$ where $1_{\cal C}$
is the unit object in the monoidal category ${\cal C}$. Then $\Delta$ can be
made into a chain complex of 
the CW complex $[0,+\infty]=\{0\}\cup (0,+\infty)\cup \{+\infty\}$.

We see that as a graded space  $\bar P_n$ is given by the formula

$$
\bar P_n=\oplus_{\lbrack T \rbrack, T\in Tree(n)}(\otimes_{v \in V_i(T)}
R_{N^{-1}(v)}\otimes \Delta^{\otimes E_i(T)})_{AutT}
$$

Since we have here a tensor product of complexes, we get the corresponding
differential $d_{\bar P}$ in $\bar P$.

\begin{prp}
The ideal $I_P$ is preserved by $d_{\bar P}$.

\end{prp}
$Proof$. Straitforward computation. $\blacksquare$

Therefore $P=P_R$ is a dg-operad which is free as a graded operad.

There is a natural morphism of dg-operads $\phi:P\rightarrow R$. In terms of
the Description 2 it can be defined as follows:

$\phi$ sends to zero all  generators of $P$ corresponding
 to trees with at least
one finite edge.
Let $T\in P$ be a tree with all infinite edges. Then $T$ gives rise
to a natural rule of composing in $R$ elements of $R_{N^{-1}(v)}$
assigned to the vertices of $T$. We define $\phi(T)\in R$ as the
result of this composition.

It is easy to check that $\phi$ is a well-defined morphism of dg-operad.

\begin{prp}
The morphism $\phi$ is a quasi-isomorphism of dg-operads.
\end{prp}

$Proof$. Follows from the spectral sequence arising from the natural
stratification of ${\cal T}$. To say it differently,
 let us consider the tautological
embedding $\psi$ of $R$ into $P$. Then $\psi$ is a right inverse
to $\phi$. It gives a splitting of $P$ into the sum $P= \psi (R)\oplus
P^{(0)}$. Here $P^{(0)}$ is spanned by the operations corresponding to trees
with finite edges only. Such a tree can be contracted to a point which
means that $P^{(0)}$ is contractible as a complex. Hence $\phi$
defines a quasi-isomorphism of complexes and dg-operads. $\blacksquare$

\subsection{Example}

Let us discuss an example when $R$ is the  operad of 
associative algebras without
the unit. We denote it by $As$.
 Then for any $n \ge 1$ we have: $As_n$ is isomorphic 
 to the  regular representation
of the symmetric group $S_n$. 

In this case the complex $P_n$ from the previous subsection can be identified
with the chain complex  of the CW-complex $K_n, n\ge 2$ described below.

The cells of $K_n$ are parametrized by $planar$ $trees$ with an additional
structure on edges.
By a 
planar tree here 
we understand a numbered tree $T$ such that for any $v\in V_i(T)$
the cardinality of $N^{-1}(v)$ is at least $2$ and this set is completely 
ordered. 
The additional structure 
is a map $E_i(T)\rightarrow \{finite,infinite\}$. We call
an edge finite or infinite according to its image under this map.
Dimension of the cell is equal to the number of finite edges of the
corresponding  planar tree.

We can either contract a finite edge or make it infinite. This defines
an incidence relation on the set of cells.

 We can picture planar trees as follows
 
 \vskip 1cm
\centerline{\epsfbox{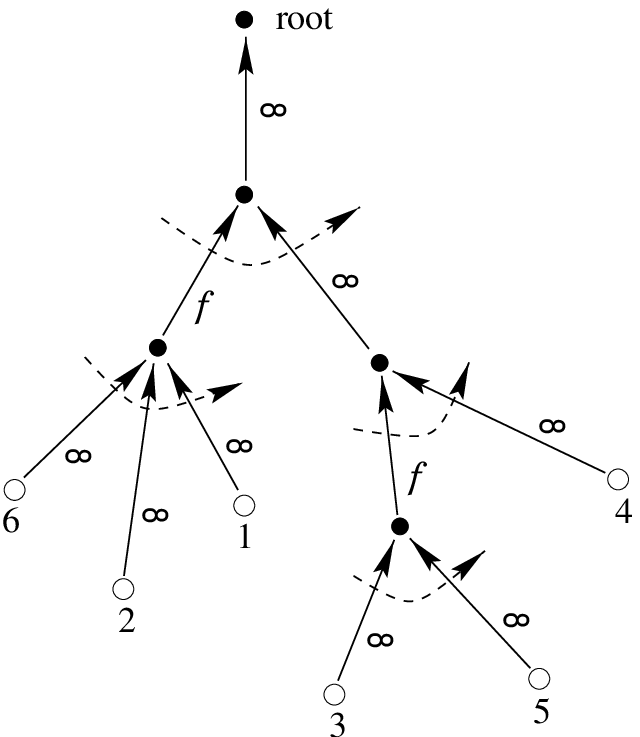}}
\vskip 1cm

 Here the dashed lines show the complete orders on set of incoming edges.
 We will not show them on other figures in the text. Instead, we will
 tacitly assume that for a given vertex the incoming edges are completely
 ordered from the left to the right.

 In this way we obtain cubical subdivisons of the Stasheff polyhedra.
 
 We depict the case $n=4$ below
 
 \vskip 1cm
\centerline{\epsfbox{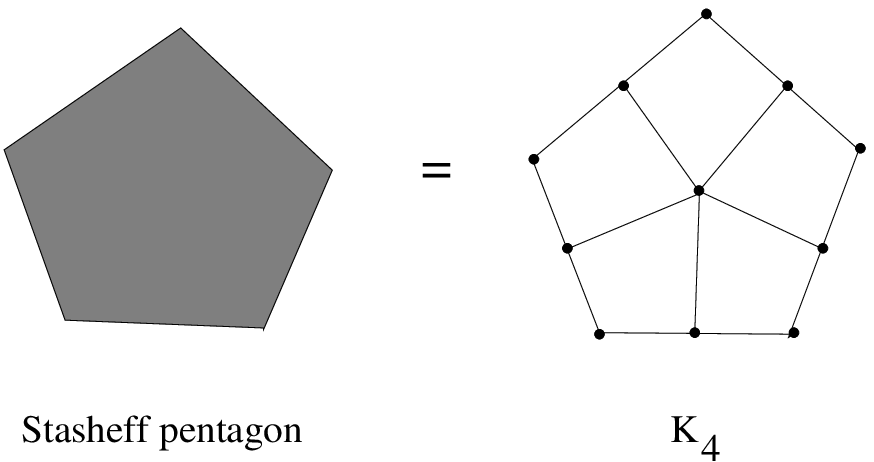}}
\vskip 1cm

\section{Minimal operad }

In this subsection we describe a dg-operad $M=(M_n)_{n\ge 1}$ which acts
naturally on the Hochschild complex of an $A_{\infty}$-algebra. We call
it $minimal$ operad. 
Let us describe this operad informally. We treat elements
$\gamma\in C^{\cdot}(A,A)$ as polylinear operations on $A$.
For given operations $\gamma_1,...,\gamma_n \in C^{\cdot}(A,A) $, and 
an $A_{\infty}$-structure $m\in  C^{\cdot}(A,A)$ we can make 
compostions in all possible ways, reading the arguments from left to right.
For example we can make an expression like this:

$\gamma(a_1\otimes...\otimes a_5)=
\gamma_2(a_1\otimes a_2\otimes m_2(a_3\otimes \gamma_1(a_4))
\otimes \gamma_3(a_5))$, etc.

Such compostions can be depicted by planar trees. The operad $M$
is spanned by operations corresponding to certain trees, which we
call admissible (see below).

It seems that the operad $M$ is close to what is described in [MS]
as ``natural transformations $(C^{\cdot}(A,A))^{\otimes n}\to C^{\cdot}(A,A)$"
(this terminology is confusing because 
the assignement $A\mapsto C^{\cdot}(A,A)$ is not a functor).

In all previous works dealing with Deligne's conjecture the authors
used operads which act on $C^{\cdot}(A,A)$ and generated by the operations
called braces. It seems that the braces generate the operad $M$, but it
is not clear what is the complete list of relations.
The advantage of our operad $M$ is that it is described directly, not
as a quotient of a free operad.

The defining properties of a dg-operad (i.e. associativity
of the composition, Leibniz rule for the differential) will become clear later.
We will give two definitions of the operad $M$. The first one is
suitable for the pure algebraic descriptions of the operadic composition
and the differential. But the signs in the formulas are not very transparent.
Second description takes care about ``parity of edges", so the 
correct signs come out automatically. In some formulas related to the
first description we will write $\pm$ having in mind that the correct
sign follows from the second description.

\subsection{Basis of $M$}

\begin{dfn} For a finite set $I$ we define an $I$-labeled planar tree 
as a triple $(T,lab, ord)$ where $T$ is a tree in the sense of Section 1.2,
$lab: I \hookrightarrow V_i(T)$ is an embedding, and $ord$ is a complete
order on the sets $N^{-1}(v), v\in V_i(T)$.

\end{dfn}

We call {\it labeled} a vertex from the image of
the map $lab$. 
All other internal vertices  are called {\it non-labeled}.
 
\begin{dfn} We will call an $I$-labeled tree admissible if it has
no tail vertices, and for 
every non-labeled internal vertex $v$ we have $|v|\ge 2$.

\end{dfn}

We denote the set of isomorphism classes of $I$-labeled planar trees by 
$Tree^{(p)}(I)$, where the upper script $p$ stays for ``planar".
 Notice that the automorphism group
of an $I$-labeled planar tree is trivial.   
For $I=\{1,...,n\}$ we will use the notation $Tree^{(p)}(n)$.

We are going to use admissible
trees unless we say otherwise.
If it will not lead to a confusion, we will simply call them trees.
Slightly abusing the notation we will 
denote an $I$-labeled tree by $T$, skipping {\it lab} and {\it ord}.
Notice that terminology here is different from the one in Section 1.5.
In particular, we  label here   internal vertices, not tails. Since we
do not consider here numbered trees (in the terminology of Section 1.5),
this change of terminology should not lead to a confusion.

We can depict trees from $Tree^{(p)}=(Tree^{(p)}(n))_{n\ge 1}$ as follows

\vskip 1cm
\centerline {\epsfbox {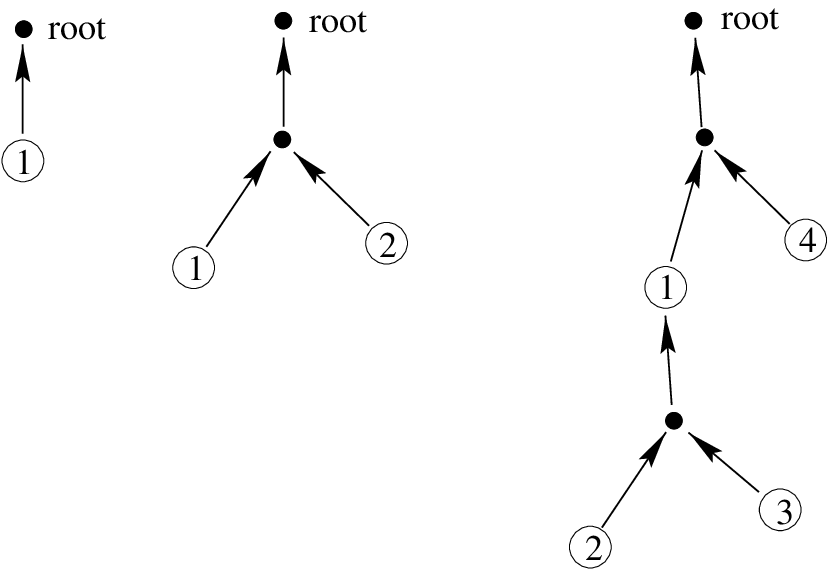} }
\vskip 1cm

Labeled vertices are depicted as circles with numbers inscribed, 
non-labeled vertices are depicted as black vertices.

We define $M_I$ to be a $k$-vector space spanned by all
elements of $Tree^{(p)}(I)$. For $I=\{1,...,n\}$
we will use the notation $M_n$. 
The symmetric group $S_n$ acts on $M_n$ permuting labeled vertices.

Abusing the notation further, we will denote
the element of $M_n$ corresponding to a tree $T$ simply by $T$.
Thus we have: $M_0=0$ and
 $M_1$ is a 1-dimensional vector space generated by $T_e$, the first tree on
the figure above. In fact $T_e$ corresponds to the unit $1_M\in M_1$
of the operad $M$.

  The  operadic
composition in $M$ and the differential will be described below.

The degree of the basis element corresponding
to a tree $T$ is equal to

$$
deg(T)=\sum_{v\in V_{lab}(T)}(-|v|)+ \sum_{v\in V_{nonl}(T)}(2-|v|)
$$
where $V_{lab}(T)$ and $V_{nonl}(T)$ denote the sets of labeled and 
non-labeled vertices respectively , and $|v|$ is the cardinality
of the set $N^{-1}(v)$. 

\subsection{Composition in $M$}

We need to
define an element of $M$ which corresponds to a tree
$T_2$ glued to a tree $T_1$ at a labeled vertex $v\in V_{lab}(T_1)$.
The trees $T_i, i=1,2$ correspond to some elements of $M$.
The resulting element will be by definition their operadic composition.
It is given
by the sum 
$$
T_1\circ_vT_2=\sum_{\beta}\pm (T_1\circ_vT_2)_{\beta}
$$
where the trees $(T_1\circ_vT_2)_{\beta}$ are defined below.

First, with the tree $T_2$ we associate a set 
$A(T_2)=\bigsqcup_{v\in V_i(T_2)}\{0,...,|v|\}$.
  We call
it the set of $angles$ of $T_2$. Obviously there is a natural map
$\kappa: A(T_2)\to V_i(T_2)$. The path in ${\bf R}^2$ which goes
from the left to the right and surrounds $T_2$ defines a complete
order on $A(T_2)$. On the following figure angles are marked by
asteriscs.

\vskip 1cm
\centerline {\epsfbox {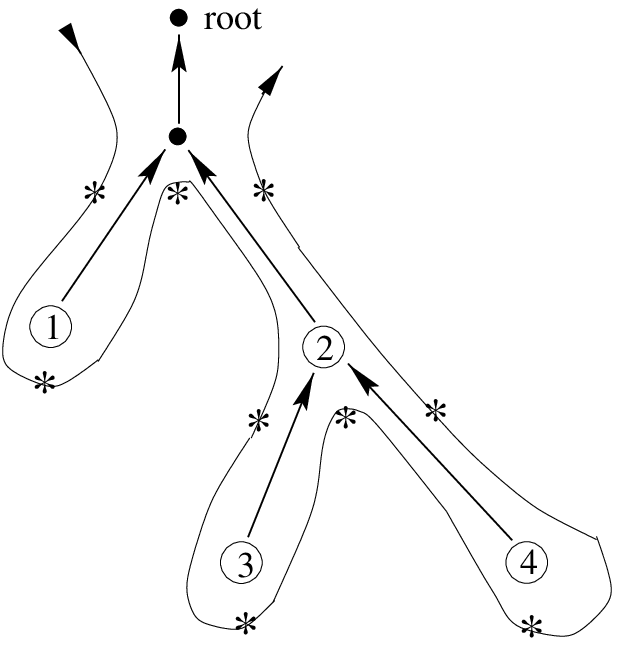} }
\vskip 1cm

 The  datum $\beta$ above is a (non-strictly)
monotonic map $\beta: N^{-1}(v) \rightarrow A(T_2)$. We will think of
this map as about the way to insert a vertex $w \in N^{-1}(v)$ and the
edge $(w,v)$ inside of an angle formed by two edges incoming to $\beta(v)$.
  
Schematically it is shown on the figure below.
  
\vskip 1cm
\centerline {\epsfbox {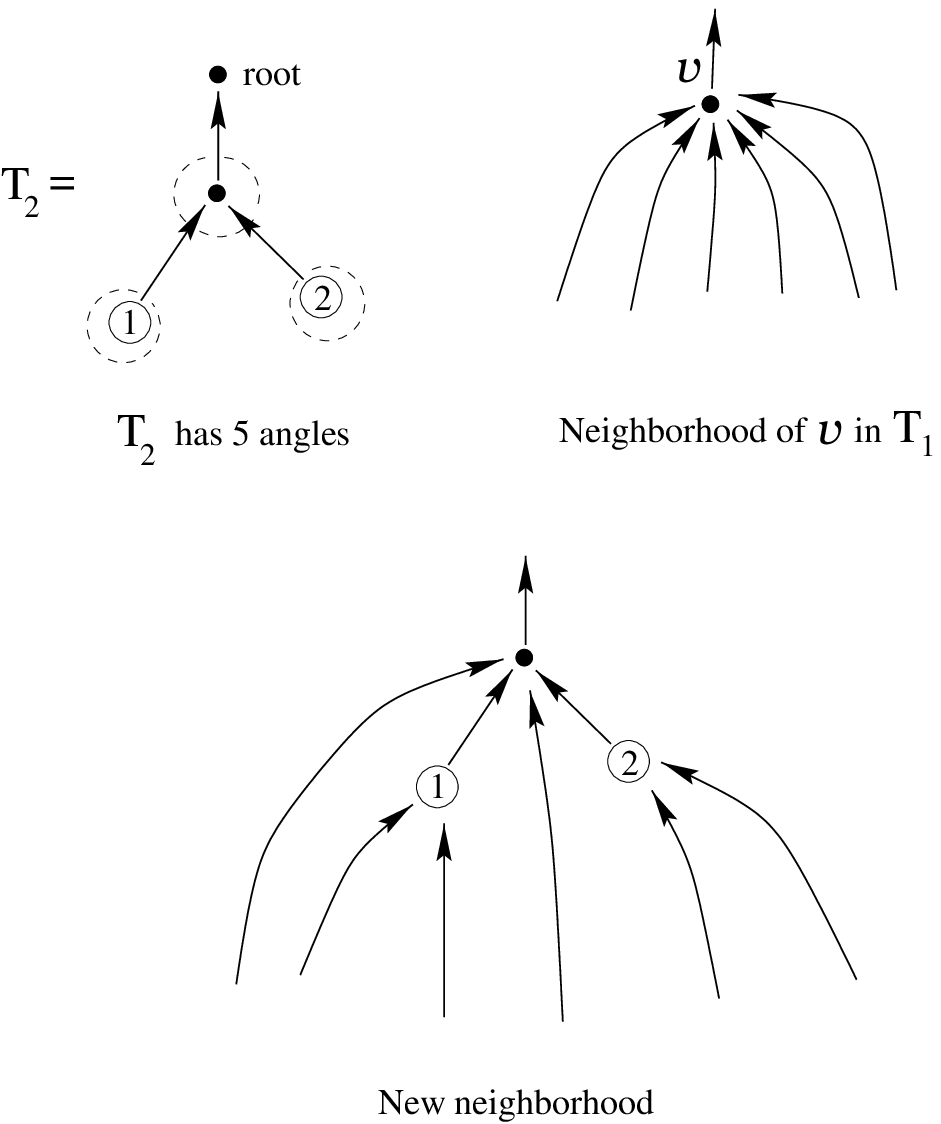} }
\vskip 1cm

Let $T_{\beta}=(T_1\circ_vT_2)_{\beta}$.
Then we define the set of vertices $V(T_{\beta})$ as $(V(T_1)\setminus\{v\})
\bigsqcup (V(T_2)\setminus \{root\})$. The map $N=N_{T_{\beta}}$ is defined 
such as follows: for all $w\in V(T_1)\setminus \{v\}$ such that $N_{T_1}(w)
\ne v$ we put $N(w)=N_{T_1}(w)$ in the self-explained notation.
Similarly if $w\in V(T_2), N_{T_2}(w)=root$ we put $N(w)=N_{T_1}(v)$.
If $w\in V(T_2), N_{T_2}(w)\ne root_{T_2}$ then we put $N(w)=N_{T_2}(w)$. 
Let us suppose that $w\in N_{T_1}^{-1}(v)$. Then we define $N(w)$ as
$\kappa(\beta(w))$. The root vertex of $T_{\beta}$ is the same as for
$T_1$. The labeling and complete orders on the sets $N^{-1}(x)$ are defined
in the natural way. Informally speaking, $T_{\beta}$ is obtained by
removing from $T_1$ the  vertex $v$ together with
all incoming edges and vertices, and gluing $T_2$ to $v$. Then we use
the map $\beta$ in order to ``insert" removed vertices.
With such a composition we obtain the structure of a graded operad on $M$.

We depict an example of the composition in $M$ below

\vskip 1cm
\centerline {\epsfbox {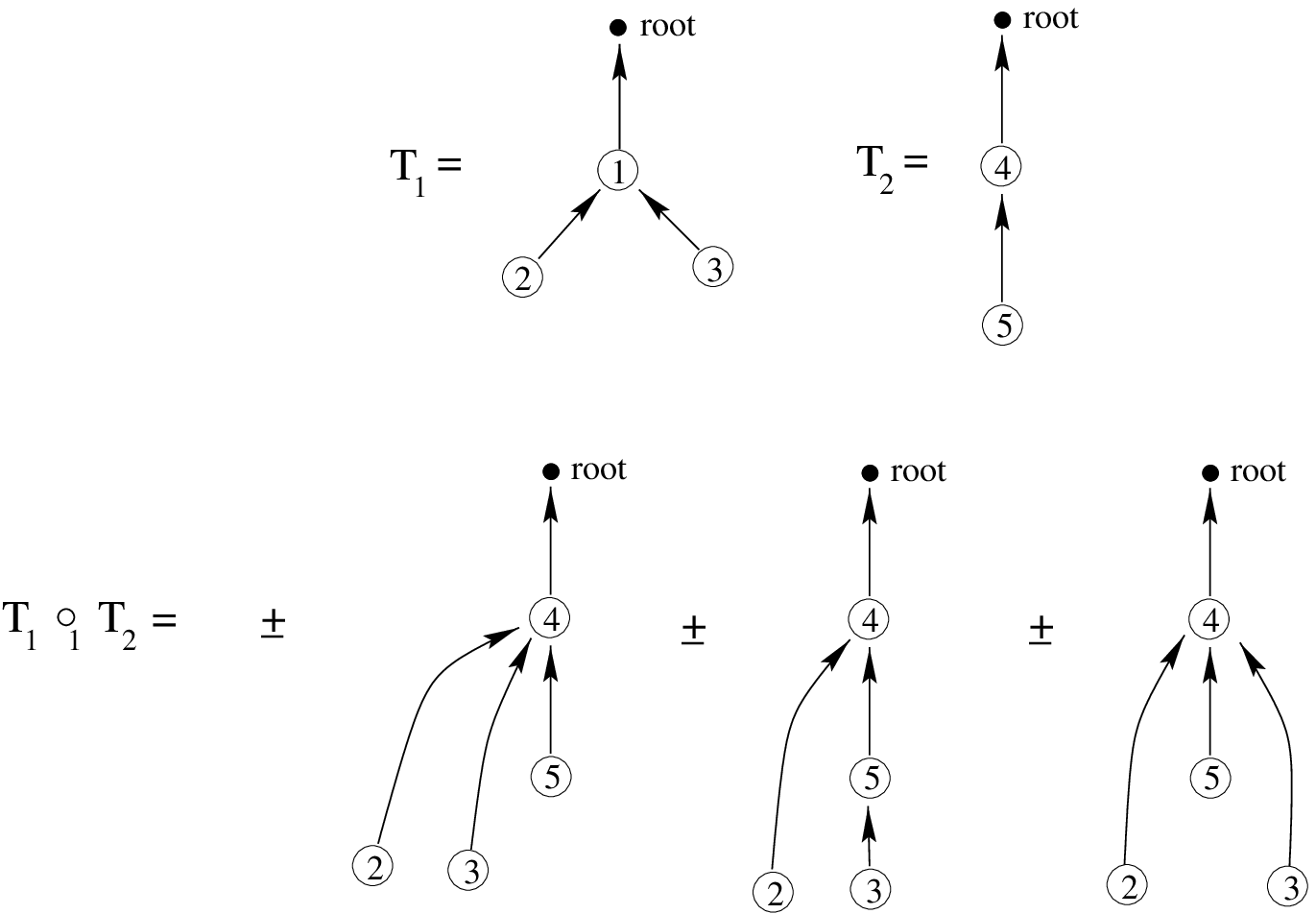} }
\vskip 1cm

\vskip 1cm
\centerline {\epsfbox {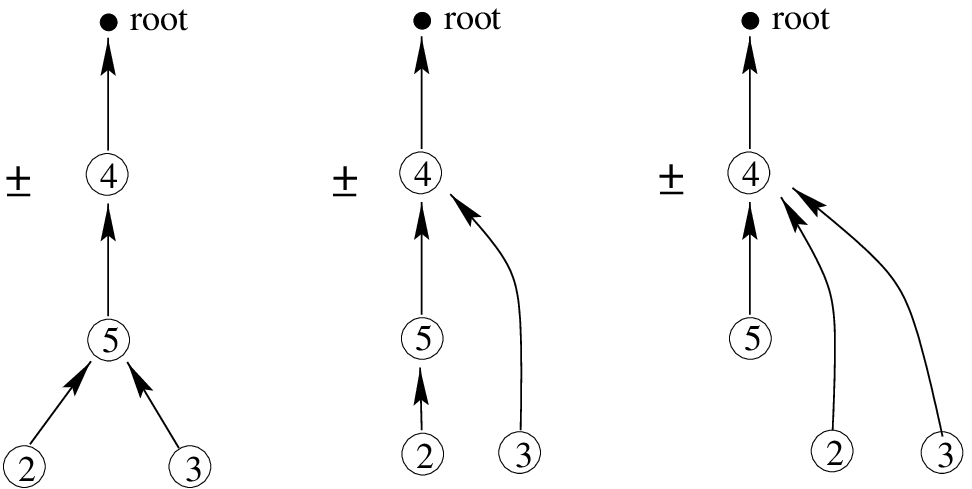} }
\vskip 1cm

\subsection{Differential in $M$}

For a generator $T$  we define
$d_M(T)=\sum_{v\in V_i(T)}d_v(T)$ where each $d_v(T)$ will be of the form
$$
d_v(T)=\sum_{i,j}\pm d_{v,i,j}(T).
$$

We need to explain the set of indices of summation and each summand.
Let us recall that for every vertex $v\in V_i(T), |v|=k$ there is a bijection
of sets $\{1,...,k\}\rightarrow N^{-1}(v)$ which defines a complete order
on the set $N^{-1}(v)$. We will identify an element of $N^{-1}(v)$ with
the corresponding number.
The indices $i,j$ in the sum above will be half-integers:
 $1/2 \le i \le j \le k+1/2$. 
 For a pair $i,j$, we define $d_{v,i,j}(T)$ to be a tree
$T^{\prime}$ such that:

a) $V(T^{\prime})=(V(T)\setminus \{v\})
\cup \{v^{up},v^{down}\}$ where $\{v^{up},v^{down}\}$ are new
vertices;

b) $root_{T^{\prime}}=root_T$;

We put 

$$
N_{T^{\prime}}(v^{down})=v^{up}, 
$$

$$
N_{T^{\prime}}(l)=v^{down},
$$
if $i<l<j$;
$$
N_{T^{\prime}}(l)=v^{up},
$$
if $l<i$ or $l>j$.

For all other vertices $u$ we put $N_{T^{\prime}}(u)=N_T(u)$.

Complete orders on the sets $N_{T^{\prime}}^{-1}(u)$ are defined in
the natural way. The tree $T^{\prime}$ has one new edge $(v^{down},v^{up})$
which we denote in pictures as ``new".

We can depict these definitions as follows

\vskip 1cm
\centerline {\epsfbox {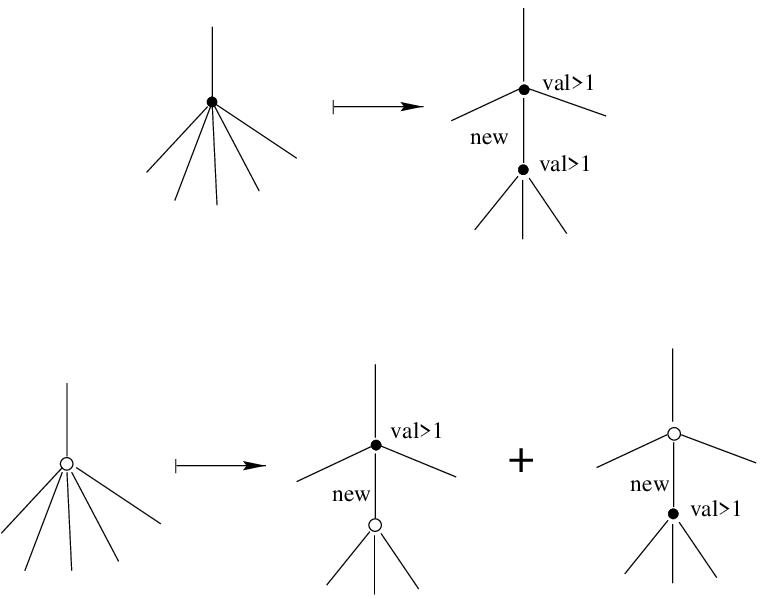} }
\vskip 1cm

The range of summation is defined differently in the case of
labeled and non-labeled vertices $v$. The idea is to keep admissible
graphs only.

1) If $v$ is non-labeled then we take  sum over $1/2\le i\le j\le k+1/2$
such that $2\le j-i$ and $1\le k-(j-i)$, in $d_{v,i,j}(T)$ both
vertices $v^{up}$ and $v^{down}$ are non-labeled.

2) If $v$ is labeled, then $d_v(T)=d_v^{(1)}(T)+d_v^{(2)}(T)$ where
the first summand $d_v^{(1)}(T)$ is the sum of 
$d_{v,i,j}(T)$ with $1\le k-(j-i)$ and
with the vertex $v^{down}$ appearing with the old label of $v$, and $v^{up}$
being non-labeled.
The second summand
$d_v^{(2)}(T)$ is the sum of $d_{v,i,j}(T)$ with $2\le j-i$, the vertex
$v^{down}$ is non-labeled, and $v^{up}$ has the same label as $v$ in $T$.
Labeling of all vertices different from $v^{up}$ and $v^{down}$
 remains the same.

\subsection{Action of $M$ on the Hochschild complex}

Let $(A,m)$ be an $A_{\infty}$-algebra.
Our next step is to define an action of a dg-operad $M$ on $C^{\cdot}(A,A)$
so the latter becomes a dg-algebra over $M$.

Let $T\in M_n, \gamma_i\in C^{\cdot}(A,A), 1\le i\le n$. We need to define
an element $T(\gamma_1,...,\gamma_n)\in C^{\cdot}(A,A)$.
It can be expressed as
a sequence of morphisms of graded vector spaces $T(\gamma_1,...,\gamma_n)_N:
A^{\otimes k}\rightarrow A$, $k=0,1,...$ 

Let $T_{(k)}$ be a unique (non-admissible) planar tree with only one
internal  vertex $v$, which is labeled, and $k$ tails, numbered
from the left to the right.

Then we have the composition

$$T_{(k)}\circ_vT=\sum_{\beta} \pm (T_{(k)}\circ_v T)_{\beta}$$
defined
in the same way as for admissible trees. Here as before $\beta:\{1,...,k\}\to
A(T_{(k)}\circ_v T)$ is a monotonic map.

Every tree $(T_{(k)}\circ_v T)_{\beta}$ has $k$ tails.
For every $v\in V_i((T_{(k)}\circ_v T)_{\beta})$ we define
a polylinear map $\gamma_v:A^{\otimes |v|}\to A$ in the following way:

1) if $v$ is labeled by $j$, $1\le j\le n$, we define $\gamma_v$
as the component of $\gamma_j$ which belongs to 
$Hom_{Vect_{\bf Z}}(A^{\otimes |v|},A)$;

2) if $v$ is non-labeled we define $\gamma_v=m_{|v|}$.

The tree $ (T_{(k)}\circ_v T)_{\beta}$ defines the way to compose 
operations $\gamma_v$ into an operation $\gamma_{\beta}:A^{\otimes k}\to A$.
We define $T(\gamma_1,...,\gamma_n)$ to be equal to the sum 
$\sum_{\beta}\pm \gamma_{\beta}$.

We claim that
in this way we get on
$C^{\cdot}(A,A)$ a structure of a dg-algebra over $M$.

This can be checked by a straightforward computation. 
It is more or less clear from the definitions
that $C^{\cdot}(A,A)$ is an algebra over the graded operad $M$. Hence the
question is about the compatibility with differentials.
The latter follows from a more general result (the Theorem $1$), 
and will be proved in 
Section 6.

\subsection{ Signs in the minimal operad}

Now we would like to discuss  the second description of the operad $M$.
It is based on the following result from the theory of
Strebel's differentials (see [St]). 

\begin{thm} Let $I=\{z_1,...,z_n\}$ be a
finite non-empty subset of the complex line ${\bf C}$.
Then there exists a unique quadratic differential
$\alpha=f(z)(dz)^2$ with
$f(z)=\sum_i {a_i\over {z-z_i}}$,
$a_i\in {\bf C}, i=1,...,n$ , $\sum_ia_i=1$,
satisfying the following property:

There exists
a unique tree $T\in Tree^{(p)}(I)$ and an embedding $j:T\to {\bf CP}^1$,
$\{root\}\mapsto \{+i\infty\}$ such that the pair
 $({\bf CP}^1\setminus j(T),\alpha)$ is equivalent (as a complex curve
 with quadratic differential) to the lower half-plane $Im(z)<0$
 equipped with the quadratic differential $i(dz)^2$.

\end{thm}

We will call such $\alpha=f(z)(dz)^2$ 
the Strebel differential associated
with the finite set $I\subset {\bf C}$.

To every Strebel differential $\alpha$ we can assign canonically
the $\{1,...,n\}$-labeled planar tree $T$. 
Conversely, having a tree $T\in Tree^{(p)}(n)$ we can look for
 sequences of pairwise distinct points $(z_1,...,z_n)$ in ${\bf C}$ 
which can appear as the sets of possible poles of the Strebel differential
from the Theorem. We remark that the set of all possible
$(z_1,...,z_n)\in {\bf C}^n\setminus \{diag\}$ form an open cell
$Str_T$. This follows from the fact that this space is a 
vector bundle of rank $2$
over a cell. The latter cell is defined by the lengths of internal edges
(lengths are taken with respect to the metric $|\alpha|$ which is well-defined
on ${\bf C}\setminus I$).
Then the embedding of $T$ is fixed up to parallel traslations. The latter span
the fiber of a vector bundle. The total space of this bundle
 is the cell $Str_T\subset {\bf C}^n\setminus \{diag\}$.
It is easy to see that $codim(Str_T)=-deg(T)=-(|E_i(T)|+2-2|V_{lab}(T)|)$.

One can show that there is a finite CW-complex 
${\Sigma}_{(n)}\subset {\bf C}^n\setminus \{diag\}$ with the cells $\Sigma_T$
labeled by $T\in Tree^{(p)}(n)$, and such that $\Sigma_T$
intersects $Str_T$ transversally and at exactly one point.

Now we can give the second definition of the operad $M$:

{\it For $n\ge 2$ the complex $M_n$ 
 is defined as the chain complex of ${\Sigma}_{(n)}$}.

Using our agreement about degrees of homological complexes,
one can check again that every complex $M_n$ is a finite-dimensional
complex concentrated in degrees $\{-(n-1),...,0\}$.

We use the second definition of  in order to derive
the following combinatorial description of the operad $M$.

One  associates to a  tree $T\in Tree^{(p)}(n)$
a one-dimensional vector space 
(in fact an abelian group, so everything can be done
over integers). It is defined by the formula
$U_T=((\widetilde {H}_{\cdot}({\bf R}))^{\ast})^{E_i(T)}\otimes 
((\widetilde {H}_{\cdot}({\bf R}^2)^{\ast})\otimes
 ((\widetilde {H}_{\cdot}({\bf R}^2))^{V_{lab}(T)}$.

In the formula $\widetilde {H}_{\cdot}(X)$ denotes the reduced homology of
the one-point compactification of $X$ (Borel-Moore homology).

The graded vector space $M_n$ coincides with the
 sum of $U_T$ over all
trees belonging to $Tree^{(p)}(n)$.

Next step is to interpret various shifts in complexes as 
tensor products with the reduced homology of vector spaces.
For example the shift by $1$ is the tensor product with the
reduced homology of ${\bf R}$.

It is more convenient to use a different notation for the same spaces
of reduced homology.

Let $L_1:=\widetilde {H}_{\cdot}({\bf R})$ and
$L_2:=\widetilde {H}_{\cdot}({\bf R}^2)$. These are $1$-dimensional
graded vector spaces of pure degrees $-1$ and $-2$ correspondingly.
For a finite set $I$  we define
a vector space

$M_I:=\oplus_{T\in Tree^{(p)}(I)}(L_1^{\ast})^{E_i(T)}\otimes L_2^{\ast}
\otimes L_2^{\otimes I}=\oplus_{T\in Tree^{(p)}(I)} U_T.$

Geometrically the new notation is related to the picture
with the Strebel differentials. We think about ${\bf R}^2$
as about direct sum of two lines ${\bf R}^2={\bf R}_{hor}
\oplus {\bf R}_{vert}$ (horizontal and vertical lines).
The vertical line is the $y$-axis in ${\bf C}={\bf R}^2$. 
 The vertical line corresponds to $L_1$
in the notation above, and the horizontal line corresponds
to $L_1^{\ast}\otimes L_2$ (we think of it as about
 ${\bf R}^2/{\bf R}_{vert}$).

There is a natural structure of a complex on $M_I$. 
We define a linear map $d_I: M_I\to M_I\otimes L_1$ as a sum

$$d_I=\sum_{T_1\rightarrow T}d_{T_1,T}.$$

Here $T_1\rightarrow T$ means that $T_1$ is obtained from $T$ by
adding an internal edge. The summand $d_{T_1,T}$ is a linear
map from $U_T$ to $U_{T_1}\otimes L_1=(U_T \otimes L_1^{\ast})\otimes L_1$.
It is given by $id_{U_T}\otimes \varepsilon_{L_1}$ where
$\varepsilon_{L_1}: {\bf 1}\to L_1^{\ast}\otimes L_1$ is the canonical
morphism in the symmetric monoidal category $Vect_{\bf Z}$.

\begin{lmm} The linear map $d_I$ defines a differential in $M_I$.

\end{lmm}

{\it Proof}. We need to prove that $(d_I\otimes id_{L_1})\circ d_I=0$.
This can be checked directly using the fact that $L_1$
has degree $-1$ , so the commutativity constraint acts on 
$L_1\otimes L_1$ as $-\sigma$ where $\sigma$ is the permutation map.
$\blacksquare$

The compostion maps for the operad $M$ can be naturally described in 
the new notation. We remark that in order to glue a tree $T_2$ to
a tree $T_1$ at a vertex $v$ first we need to define a sequence of
trees $\{(T_1\circ_v T_2)_{\beta}\}$
 (they were also denoted by $T_{\beta}$ in the
previous description of the gluing). For every $(T_1\circ_v T_2)_{\beta}$ the 
set of internal edges is the disjoint union of the corresponding
sets for $T_1$ and $T_2$. The set of labeled vertices for
$(T_1\circ_v T_2)_{\beta}$ is $(V_{lab}(T_1)\cup V_{lab}(T_2))\setminus\{v\}$.
Then we need to define a morphism of vector spaces
$U_{T_1}\otimes U_{T_2}\to U_{(T_1\circ_v T_2)_{\beta}}$. 

Equivalently we need to define a  morphism of the graded vector space

$(L_1^{\ast})^{\otimes E_i(T_1)}\otimes L_2^{\ast}\otimes L_2^{\otimes I_1}
\otimes (L_1^{\ast})^{\otimes E_i(T_2)}\otimes L_2^{\ast}
\otimes L_2^{\otimes I_2}$
to the graded vector space 
$(L_1^{\ast})^{\otimes (E_i(T_1)\sqcup E_i(T_2))}\otimes L_2^{\ast}
\otimes L_2^{\otimes ((I_1\setminus \{v\})\sqcup I_2)}.$
(The notation is self-explained).

We define the morphism to be the identity  on the tensor factors
marked by  same edges or same elements of the sets $I_k, k=1,2$.

For example
 $(L_1^{\ast})^{\otimes E_i(T_1)}\otimes (L_1^{\ast})^{\otimes E_i(T_2)}$
is identically mapped to 
$(L_1^{\ast})^{\otimes (E_i(T_1)\sqcup E_i(T_2))}$.

Hence it is enough
to define a morphism $L_2^{\ast}\otimes L_2^{\ast}\otimes L_2
\to L_2^{\ast}$, where first and third tensor factors in
the LHS correspond to the tree $T_1$, middle tensor factor
corresponds to the tree $T_2$ and the space $L_2^{\ast}$ in
the RHS corresponds to the tree $(T_1\circ_v T_2)_{\beta}$. Then one uses
the canonical evaluation map $ev: L_2^{\ast}\otimes L_2\to {\bf 1}$
 in the symmetric monoidal category $Vect_{\bf Z}$.

\subsection{Signs in the Hochschild complex}

In this subsection we are going to reformulate the definition
of the Hochschild complex 
using the reduced homology spaces $L_1$ and $L_2$. Then one
gets the sign agreement with the computations in the previous subsection.

For a graded vector space $A$ we define

$$C=C^{\cdot}(A,A)=\oplus_I \underline{Hom}(A^{\otimes I}, A)
\otimes (L_2^{\ast}\otimes L_1)^{\otimes I}$$
where the sum is taken over all non-empty completly ordered finite sets,
and $\underline{Hom}$ is the internal $Hom$ in the 
tensor category $Vect_{\bf Z}$.

Let us demonstrate how the structures in $C$ can be reformulated
by means of this language.

The Gerstenhaber bracket is a map 
$C\otimes C\to C\otimes (L_2^{\ast}\otimes L_1)$.
Now the multiplication $m$ can be described as a point
of a dg-manifold: $m\in C\otimes L_2$. The bracket $[m,m]$
defines a point of the dg-manifold
$C\otimes L_2^{\ast}\otimes L_1\otimes L_2\otimes L_2=
(C\otimes L_2)\otimes L_1$. 
If $[m,m]=0$ then we have a differential $d_m=[m,\cdot]:C\to C\otimes L_1$.

The corresponding structure of DGLA on $C$ can be reformulated 
such as follows.

 Consider the class ${\cal F}$ of $simple$ $forests$.
 A simple forest $F$ is
 a finite collection of planar trees 
 $F=\{T_{\alpha}\}_{\alpha\in \Omega}$ with no internal
edges. For every $F\in {\cal F}$ we define the graded vector space

$$W_F:=\otimes_{\alpha\in \Omega}Hom(A^{\otimes V_t(T_{\alpha})},A)
\otimes (L_2^{\ast}\otimes L_1)^{\otimes V_t(T_{\alpha})}\otimes L_2 .$$

We have the natural groupoid structure on ${\cal F}$. The correspondence
$F\mapsto W_F$ defines a functor from this groupoid to the category
$Vect_{\bf Z}$. The free cocommutative coalgebra cogenerated
by $C\otimes L_2$ can be described as a colimit of this functor:
$Coalg(C\otimes L_2)=colim_{F\in {\cal F}}W_F=\oplus_{F/iso}(W_F)_{AutF}$.

The coalgebra structure on $Coalg(C\otimes L_2)$ can be described
in these terms.
Namely
$\Delta\circ pr=\sum_{F_1\subset F}(pr\otimes pr)\circ i(F_1,F_2)$.

Let us explain the notation.

Here $\Delta$ is the coproduct on the free cocommutative coalgebra cogenerated
by $C\otimes L_2$, $pr$ is the projection to the coinvariants.
The sum is taken over all subforests $F_1$
(unions of some connected components of $F$), and we fix the splitting
$i(F_1,F_2): W_F\simeq W_{F_1}\otimes W_{F_2}$, where $F_2$ is the
complementary forest. 
We leave to the reader straightforward reformulations
of other structures on the Hochschild
complex and  checking the signs.

\begin{rmk} For an $A_{\infty}$-algebra $A$ one can define
the opposite algebra $A^{op}$. In the geometric language
of this section it corresponds to the antipodal involution
on ${\bf R}_{hor}$.

\end{rmk}

\section{Morphism of dg-manifolds}

\subsection{Generators of the free operad}

The main purpose of this section is to prove the Theorem 1.
It will be done in subsection 6.2. This subsection is
devoted to some technical preparatory material.

Applying the general theory of Section 4 to
the case of the operad $M$ we obtain its free resolution $P$.
It is a dg-operad which is free as a graded operad. Its graded
components $P_n$ can be described explicitly in terms of
the operad $M$. We know that  $P_0=0$ and $P_1=k\cdot 1_P$.
 One can easily describe  the space 
 $G_n=(Free^{\prime}(M^{\prime}\lbrack 1\rbrack)\lbrack -1 \rbrack)_n$ 
of generators of $P_n, n\ge 2$.

Namely, $G_n$ is a direct sum of $1$-dimensional graded vector spaces
$W_{\bf T}$.
The sum is taken over the set of equivalence classes of collections
${\bf T}=(T,\{T_v\}_{v\in V_i(T)})$ where:

a) $T\in Tree(n), n\ge 2$ such that for every $v\in V_i(T)$ we
have $|v|\ge 2$;

b) $T_v \in Tree^{(p)}(N_T^{-1}(v))$.

In order to be consistent with the notation of Section 7,
we should also label all internal edges of ${\bf T}$ by 
an additional label {\it finite}.
 We omit the labeling since it will not
be used in this section.

Clearly the automorphism group of any such ${\bf T}$ is trivial.
The symmetric group $S_n$ acts freely on $G_n$
permuting the tails of $T$. Any ${\bf T}$ gives rise to
the $1$-dimensional vector space 
$W_{\bf T}=L_1^{\otimes E_i(T)}\otimes U_{T_v}$,
where the space $U_{T^{\prime}}$ was defined in 5.5.
In what follows we will identify ${\bf T}$ with the
corresponding (up to a sign) generator of $G_n$.

Notice that if $E_i(T)=\emptyset$ then $T$ has the only 
internal vertex $v$ and $W_{\bf T}=U_{T_v}$. Moreover, the set
$N^{-1}(v)$ is naturally identified with the set $\{1,...,n\}$, so
that we can write $T_v\in Tree^{(p)}(n)$.

There is a natural morphism of $S_n$-modules $pr_n:G_n\to M_n$,
such that $pr_n(W_{\bf T})=0$ if $E_i(T)=\emptyset$, otherwise
$pr_n(W_{\bf T})=id_{T_v}$ where $v$ is the only internal
vertex of $T$.

A generator ${\bf T}$ is pictured  by a tree  $T$ with $n$ numbered tails,
and with generators of
 $M$  inscribed into all internal vertices of $T$. If $T_v\in M$
is inscribed into a vertex $v\in V_i(T)$ then the cardinality
of $V_{lab}(T_v)$ is equal to the cardinality of $N^{-1}(v)$.

Now we can return to the Theorem 1. We have constructed the minimal 
operad $M$. The Hochschild complex $C^{\cdot}(A,A)$ of
an $A_{\infty}$-algebra $A$ is an algebra over $M$
 (if we forget the differentials).
 Then the free resolution $P$ of $M$ acts on
$C^{\cdot}(A,A)$,  so the latter becomes an algebra over 
the graded operad $P$ (again we forget the differentials).

On the other hand, let $A$ be a graded vector space and
 $m\in C=C^{\cdot}(A,A)$,
but not necessarily $\lbrack m,m \rbrack=0$. Then the constructions of the
 Section 5 give rise to a sequence of elements 
$\rho (m)=(\rho (m)_n)_{n\ge 1}$ of 
$\underline{Hom}
(M_n\otimes_{S_n}C^{\otimes n},C)$ (no conditions on $m$ were
imposed by the construction). The following lemma is easy to prove

\begin{lmm} The sequence $\rho (m)$ defines a structure of 
a graded $M$-operad on $C$.

\end{lmm}

Compositions $\gamma(m)_n=\rho (m)_n\circ pr_n$
for $n\ge 2$ give rise to a sequence of
$S_n$-equivariant maps $G_n\otimes_{S_n}C\to C$. Since
the graded operad $P$ is freely generated by $G=(G_n)_{n\ge 2}$,
these compositions define a point in ${\cal M}(P,C)$.
We denote the sequence $(\gamma(m)_n)_{n\ge 2}$ by $\gamma(m)$.
Thus we can write $\gamma(m)=\rho (m)\cdot pr$ where $pr=(pr_n)_{n\ge 2}$.

We can  define an element $d_m\in
\underline{Hom}(C,C)\lbrack 1\rbrack$ in the natural way:
 $d_m=\lbrack m,\cdot\rbrack$. In general $d_m^2\ne 0$.
 In any case
it can be naturally extended to a map 
$f:{\cal M}_{cat}({\cal A}_{\infty},A)\rightarrow
{\cal M}(P,C)$ such that $f(m)=(d_m,\gamma(m))$.

We claim that this map is a morphism of dg-manifolds
and it satisfies the condition $p\cdot f=\nu$ of the Theorem 1. 
It will be proved in the next subsection.

\subsection{Proof of the theorem}

 In the course of the proof  we will
not pay much attention to the signs. The reason for that was
explained in the previous section. Namely, our second
description of the operad $M$ (with reduced homology)
 gives automatically the agreement of signs.

 Proof of the theorem will occupy the rest of this subsection. We start 
 with some general considerations.

 Let us recall that points of the moduli space $Y={\cal M}(P,C)$ parametrize
 pairs $(d_C,\gamma)$ where
 $d_C: C\rightarrow C\otimes L_1$ is a morphism in $Vect_{\bf Z}$, 
 and $\gamma$ is an action of the space of generators 
 $G=(G_n)_{n\ge 2}$ on $C$. Having the action $\gamma(m)=\rho (m)\cdot pr$
 of $P$ on $C$ we would like to compute the odd vector field $d_Y$
 at the point $(d_m,\gamma (m))$. Clearly 
 $d_Y(d_C,\gamma)=(d_C^2,\bar \gamma)$ where 
 $\bar \gamma=(\bar \gamma_2,
 \bar \gamma_3,...)$ and 
 $\bar \gamma_n\in \underline{Hom}(G_n\otimes_{S_n}C^{\otimes n},C\otimes L_1)$.
 
 We need to compute $\bar \gamma$. One can  extend $d_C$
 to $\underline{Hom}(C^{\otimes n},C)$ using the Leibniz rule.
 We denote this extension by
  $d_C^{(n)}$. To write down $\bar \gamma$ we need to know  the images
 of the generators of $P$ under all $\gamma_n$.
  We know that the maps $pr_n$ send
 to zero all of them except of the trees with one non-labeled vertex and
 $n$ numbered tail vertices. The latter is depicted below
 
\vskip 1cm
\centerline {\epsfbox {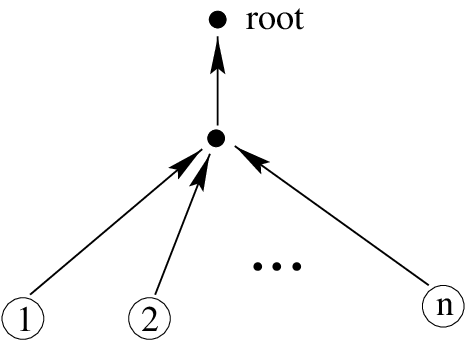} }
\vskip 1cm

 Then the direct computation together the previous lemma
 show that the following result holds
 
 \begin{lmm}
 The component $\bar{\gamma}_n$ of $d_Y(f(m))$ is equal to
  
  $$(d^{(n)}_m\cdot \rho(m)_n+\rho(m)_n\cdot d_{M_n})\cdot pr_n$$ 
 
 where $d_{M_n}$ is the differential
 $d_M$ being restricted to $M_n$.
 
 \end{lmm}
{\it Proof}. We sketch the proof. First we observe that all components
of $\bar{\gamma}_n$, which corresponds to the generators of $G_n$
with more than two internal edges, vanish. Using the previous lemma
one can show that the generators with one internal edge give no input 
as well. After that one can make computations with elements of the
operad $M$ only. Then the direct computation proves the lemma. $\blacksquare$

  We would like to prove that
 the map $f:X={\cal M}_{cat}({\cal A}_{\infty},A) \rightarrow Y={\cal M}(P,C)$
  is a dg-map, i.e.
 it transforms the odd vector field $d_X$ into the odd vector field 
 $d_Y$. Vector field on $X$ is given by $\dot m=d_X(m)={1\over2}
 \lbrack m,m\rbrack$.

The image of the map $f$ belongs to the vector subspace
$\underline{Hom}(C,C)[1]\oplus \underline{Hom}(M^{\prime}(C),C)$ of $Y$.
Therefore the image $\dot f(m)=f_{\ast}(d_X(m))$ of the tangent
vector $d_X(m)$ belongs to the same vector space. Thus, we have to show
that the second component $\dot{\gamma}$ of $\dot f(m)$ is equal to 
$\overline {\gamma}$. 

 Its first component is equal to
 ${1\over2}\lbrack \dot m,\cdot \rbrack$ which is
 the same as ${1\over2}\lbrack\lbrack m,m\rbrack ,\cdot\rbrack=d_m^2$.
 Notice that $\overline {\gamma}-\dot{\gamma}$
 can be
 considered as  an action of an action of
 $M^{\prime}$ on $C$. We decompose it into
 the sum of terms corresponding to planar trees from $Tree^{(p)}$.
 For such a tree $T$ the component of
 $\overline {\gamma}-\dot{\gamma}$ is a sum of four terms
 described below.
 
 A. These terms correspond to the differential $d_M$ of the operad $M$. We can
 schematically write them as $(d_M T)((c_i), m)$.
  We are inserting $c_i \in C$ in labeled
 vertices and $m$ in non-labeled vertices.
 
 B. These terms can be schematically written as 
 $\sum_i T(\lbrack m,c_i\rbrack,
 (c_j)_{j\ne i}, m)$. We are inserting $m$ in non-labeled vertices,
 elements $c_j$ and $\lbrack m,c_i\rbrack$ in the corresponding labeled
 vertices.
 
 C. These terms can be schematically written as $\lbrack m,T((c_i),
 m) \rbrack$ in the notation above. These terms appear when we apply
 the differential (= commutator with $m$) to the tree with $c_i$ inserted
 in labeled vertices and $m$ inserted in non-labeled vertices.

D. These terms correspond to $-\dot{\gamma}$.
Each of them consists of the replacement of $m$ in one
 non-labeled vertex $v$ by $\dot m$.
 The latter can be in turn replaced by ${1\over2}\lbrack m,m\rbrack$. We can
 schematically write the resulting sum as 
 $\sum_{v\in V_{nonl}(T)}T((c_i),{1\over2}\lbrack m,m \rbrack)$.
 We are inserting the element $m$ in all non-labeled
 vertices of $T$ except $v$, the element
  ${1\over2}\lbrack m, m\rbrack$ in the non-labeled vertex $v$ and $c_i$ to all
  labeled vertices.
 
Notice that although a tree $T$ is always admissible, planar
trees appearing in the decomposition of  $\dot{\gamma}$ are not
necessarily admissible. This means that valencies of some 
non-labeled vertices can be either $0$ or $1$.

In order to depict all four cases we use the following notation: composition
of an operation $\alpha$ sitting in a vertex $v$ with $m$ produces a new
tree with a new non-labeled vertex $w$, as well as
a new edge $(w,v)$ where $N(w)=v$. Similarly, a composition of $m$ with $\alpha$
produces a new tree with a non-labeled vertex $w$, as
well as a new edge $(v,w)$ such that $N(v)=w$. The commutator $\lbrack \alpha,
m \rbrack$ corresponds to the difference of the above-mentioned trees.
This agreement will be used also in the case when $\alpha=m$ thus giving
the way to depict trees with ${1\over2}\lbrack m,m\rbrack$ inserted. 
In the pictures
below we show neighborhoods of 
vertices where the original tree changes.

\vskip 1cm
\centerline {\epsfbox {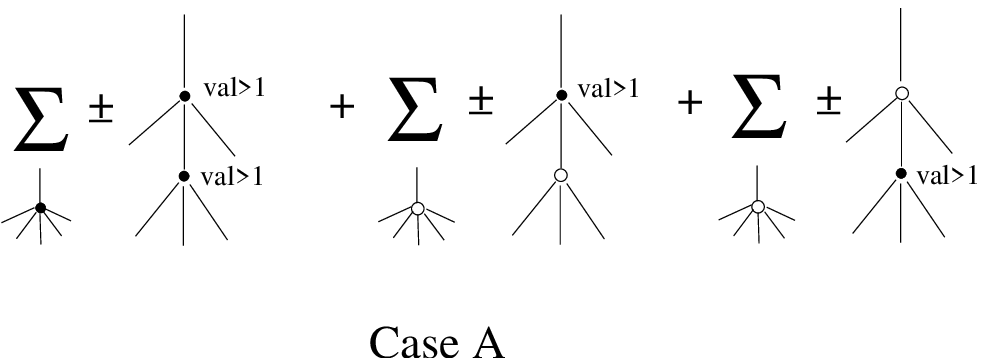} }
\vskip 1cm

\vskip 1cm
\centerline {\epsfbox {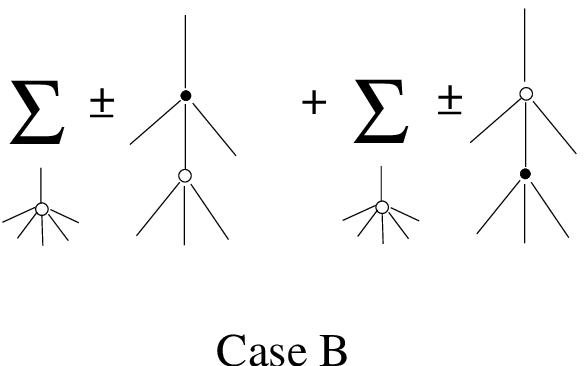} }
\vskip 1cm

\vskip 1cm
\centerline {\epsfbox {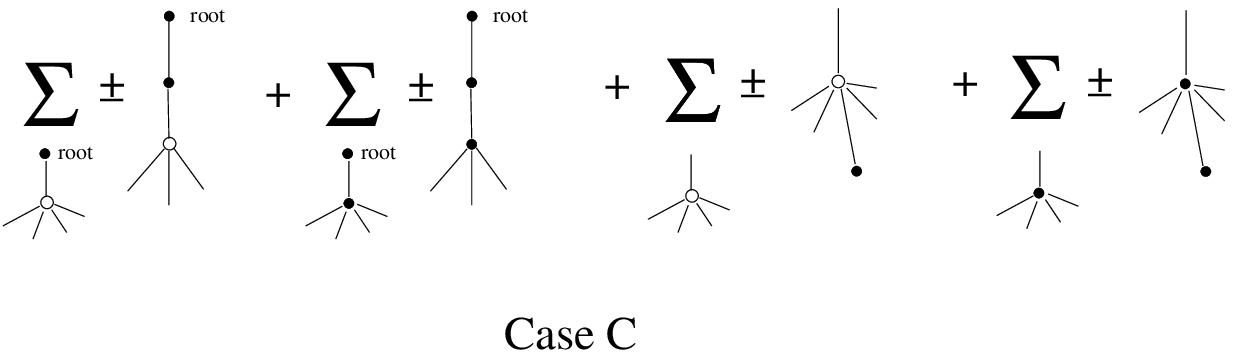} }
\vskip 1cm

\vskip 1cm
\centerline {\epsfbox {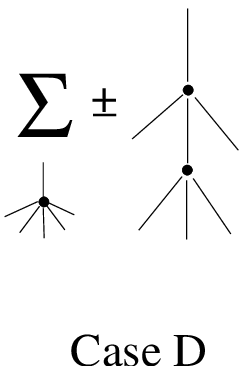} }
\vskip 1cm

Then we split these terms in the following way:

$A=A_1+A_2+A_3$, where:

a) the terms $A_1$ correspond to the two adjoint
non-labeled vertices of valency $\ge 2$ with $m$ inserted in each;

b) the terms $A_2$ correspond to the labeled vertex adjoint to a
non-labeled one of valency $\ge 2$ with $m$ inserted in the latter;

c) the terms $A_3$ correspond to the non-labeled vertex of valency $\ge 2$
with  $m$ inserted in it adjoint to a labeled vertex.

These cases can be depicted as follows

\vskip 1cm
\centerline {\epsfbox {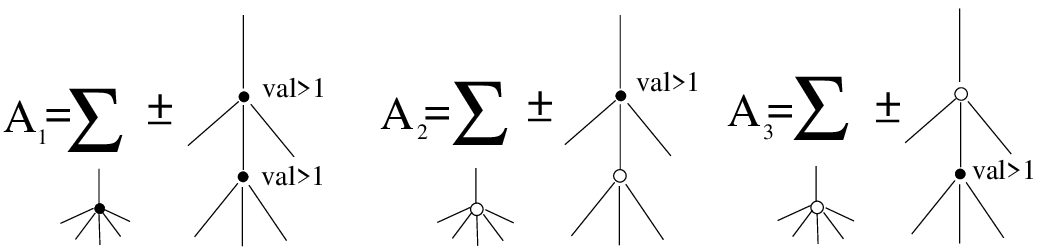} }
\vskip 1cm

Similarly we split the terms $B$ such as follows:
$B=B_{+,1}+B_{+,2}+B_{-,0}+B_{-,1}+B_{-,2}$, where individual summands
are depicted below.

\vskip 1cm
\centerline {\epsfbox {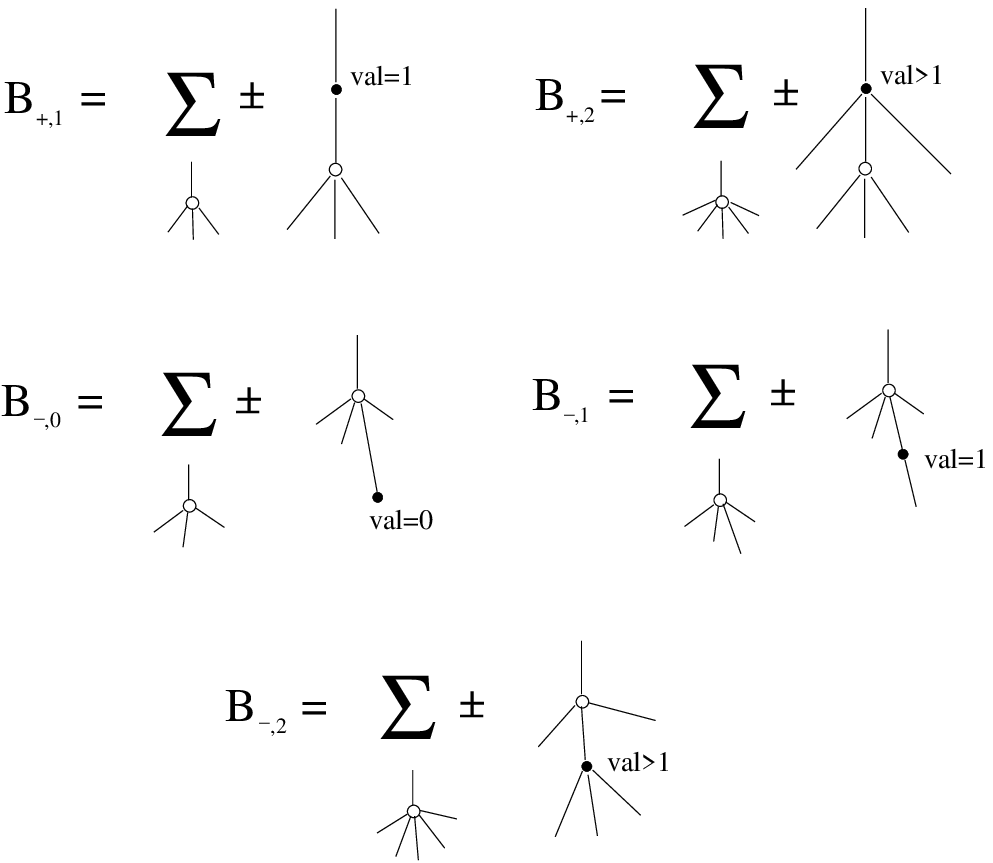} }
\vskip 1cm

We split the terms $C$ such as follows: 
$C=C_{root}^{\prime}+C_{root}^{\prime\prime}+C_{\circ}+C_{\bullet}$
where individual summands are depicted below.

\vskip 1cm
\centerline {\epsfbox {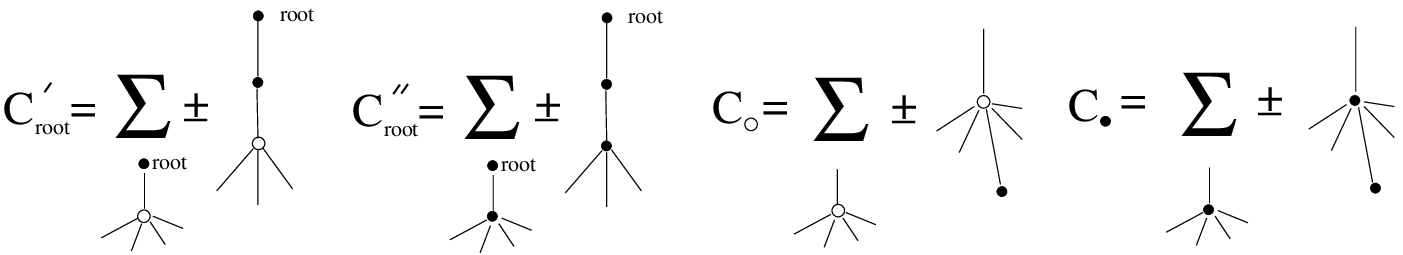} }
\vskip 1cm

We split the terms $D$ such as follows: $D=D_0+D_{1,A}+D_{1,B}+D_2$
where individual summands are depicted below.

\vskip 1cm
\centerline {\epsfbox {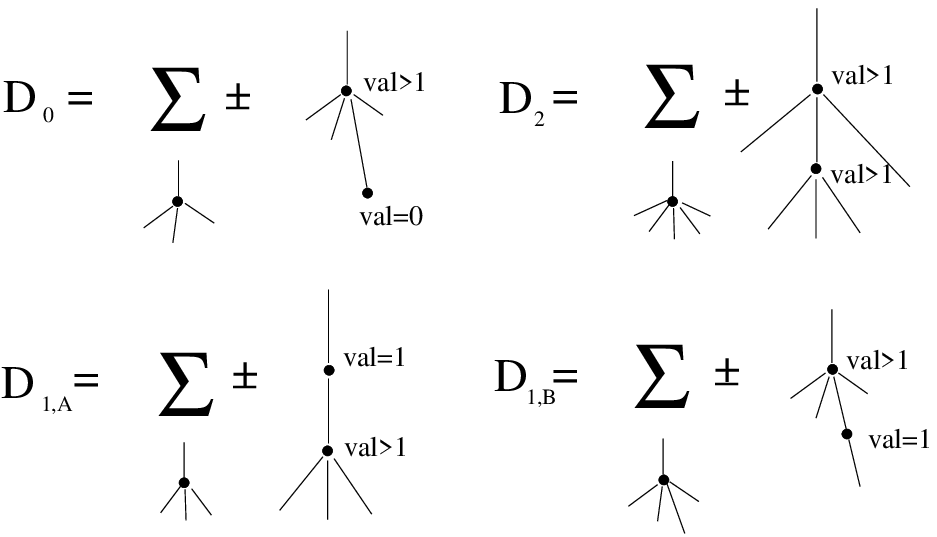} }
\vskip 1cm

We see that $A_1+D_2=0$, $A_2+B_{+,2}=0$, $A_3+B_{-,2}=0$,
 $B_{-,0}+C_{\circ}=0$, $D_0+C_{\bullet}=0$.

Furthermore we can split each of the remaining terms into summands
 $B_{+,1}=B_{+,1}^{\prime}+B_{+,1}^{\prime\prime}+B_{+,1}^{\prime\prime\prime},
  B_{-,1}=B_{-,1}^{\prime}+B_{-,1}^{\prime\prime}$, 
 $D_{1,A}=D_{1,A}^{\prime}+D_{1,A}^{\prime\prime}+
D_{1,A}^{\prime\prime\prime}$,
$D_{1,B}=D_{1,B}^{\prime}+D_{1,B}^{\prime\prime}$. 
This splitting is depicted below.

\vskip 1cm
\centerline {\epsfbox {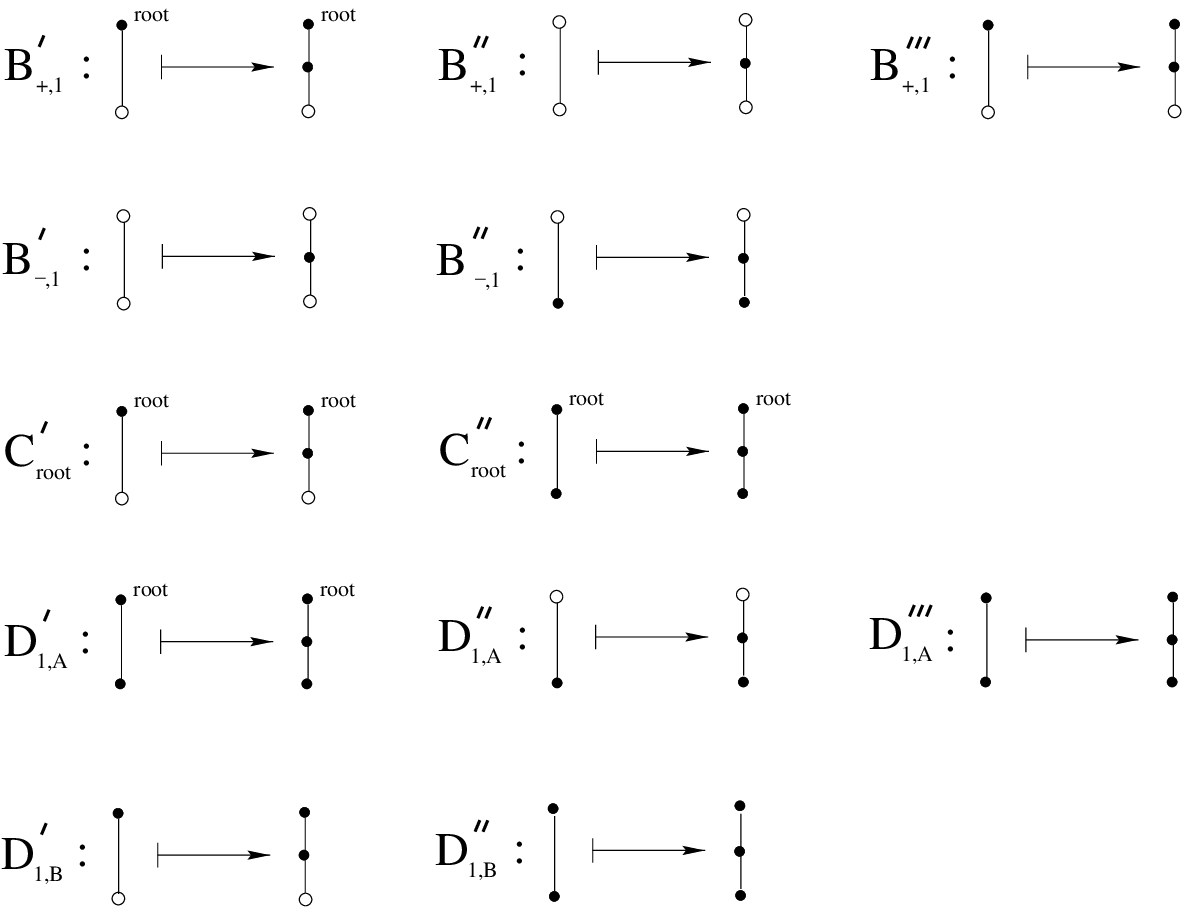} }
\vskip 1cm

We combine these summands into  six groups corresponding to the six
types of edges depicted below

\vskip 1cm
\centerline {\epsfbox {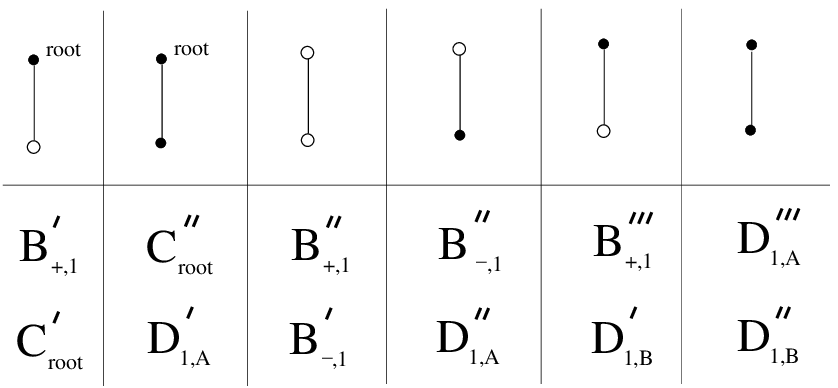} }
\vskip 1cm

We see that $B_{+,1}^{\prime}+C_{root}^{\prime}=0, C_{root}^{\prime\prime}+
D_{1,A}^{\prime}=0, B_{+,1}^{\prime\prime}+B_{-,1}^{\prime}=0$, 
$B_{-,1}^{\prime\prime}+D_{1,A}^{\prime\prime}=0, B_{+,1}^{\prime\prime\prime}
+D_{1,B}^{\prime}=0, D_{1,A}^{\prime\prime\prime}+D_{1,B}^{\prime\prime}=0$.

Then we conclude that $B_{+,1}+C_{root}^{\prime}
+C_{root}^{\prime\prime}+D_{1,A}+B_{-,1}+D_{1,B}=0$.
Thus $A+B+C+D=0$. This means that the map $f:X\to Y$ is
a dg-map. Obviously it is $GL(A)$-equivariant and
satisfies the condition $pf=\nu$.

This concludes the proof of the Theorem 1. 

\subsection{Remark about a generalization}

Let $A$ be an $A_{\infty}$-algebra, and $C=C^{\cdot}(A,A)$
be its Hochschild complex. Admissible planar trees with $n$ labeled
vertices give rise to operations
$C^{\otimes n}\to C, n\ge 1$.
Analogously, planar trees with $n$ labeled vertices and $m$ tails,
$n\ge 0, m\ge 1$ give rise to operations 
$C^{\otimes n}\otimes A^{\otimes m}\to A$.

Let us restrict ourselves to such trees that $|v|\ge 2$ 
for all
non-labeled internal vertices $v$. 
Thus we obtain a colored operad $M^{(2)}$ with two colors
{\it (Alg, Hoch)} acting on 
the set of pairs
$(A,C^{\cdot}(A,A))$, where $A$ is an $A_{\infty}$-algebra.
Clearly $M^{(2)}$ contains as suboperads both ${\cal A}_{\infty}$
and $M$.
Presumably a result analogous to the Theorem $1$ holds for the colored
operad $M^{(2)}$.

\section{Deligne's conjecture}

\subsection{Preliminaries}

We are going to prove the following result. 

\begin{thm} Let $P$ be the free  resolution of the minimal operad $M$
constructed in Section 6. Then there is a  homomorphism
of dg-operads $P\to Chains(FM_2)$ which induces
an isomorphism on cohomology (i.e. it is
a quasi-isomorphism of dg-operads). 
Here $Chains(FM_2)$ is the 
 chain operad for the 
Fulton-Macpherson operad of configurations of points in ${\bf R}^2$
(see Section 7.2 below). 

\end{thm}

We will construct a homomorphism
of the dg-operads which induces a quasi-isomorphism of the chain
complexes. Such a homomorphism is not defined canonically.
Different choices are naturally parametrized by a contractible
topological space.

Using the fact that the operad $Chains(FM_2)$ is free as an operad (not
as dg-operad), one can invert the quasi-isomorphism mentioned
in the theorem. Since the Hochschild complex is
a $P$-algebra, we obtain the following Corollary
known as Deligne's conjecture (see for example [Ko3], [V], [MS]).

\begin{cor} a) The Hochschild complex $C^{\cdot}(A,A)$
 of an $A_{\infty}$-algebra $A$
can be equipped with a structure of an algebra over the operad $Chains(FM_2)$.

b) The corresponding structure of a $H^{\cdot}(Chains(FM_2))$-algebra
on the Hochschild cohomology $H^{\cdot}(A,A)$ coincides with
the standard structure of a Gerstenhaber algebra on the Hochschild
cohomology of an $A_{\infty}$-algebra.

Same results remain true with $Chains(FM_2)$
being replaced by the operad $Chains(E_2)$ of chains on 
the little disc operad $E_2$ (see Section 7.2 for the definition). 

\end{cor}

\begin{rmk} It is easy to see that the operad $P$ 
can be defined over ${\bf Z}$.
It follows that such an operad acts on the Hochschild complex
of an $A_{\infty}$-algebra defined over a field of arbitrary characteristic.
The quasi-isomorphism $P\to Chains(FM_2)$  can be 
also defined over the ring of integers (see also [MS]).

\end{rmk}

The proof of the Theorem and the Corollary will occupy the rest
of the section.

We are going to use the following strategy.

Let ${\bf T}$ be a ``meta-tree" corresponding to a generator 
of $P$ (see Section 6 and Definition 17 below).

1) To every  ${\bf T}$ we are going to associate
a contractible closed subspace $X_{\bf T}$ 
in the Fulton-Macpherson compactification
of the configuration space of points in ${\bf R}^2$ modulo shifts
and dilations.

2) The collection of subspaces $(X)_{{\bf T}}$ will
 satisfy the following properties:

a) The correspondence ${\bf T}\mapsto X_{\bf T}$ 
is $S_n$-equivariant;

b) if a tree ${\bf T}^{\prime}$ appears as a summand in the formula for
 $d_P({\bf T})$ then $X_{{\bf T}^{\prime}}\subset X_{{\bf T}}$;

c) if a composition of trees ${\bf T}_1\circ {\bf T}_2\in P$ 
appears as a summand 
in the 
formula for $d_P( {\bf T})$ then the operadic composition 
$X_{{\bf T}_1}\circ X_{{\bf T}_2}$
belongs to the stratum $X_{\bf T}$. The latter composition of the strata
has meaning because the Fulton-Macpherson compactifications form
a topological operad.

3) For any generator ${\bf T}$ of the operad $P$
we will choose inductively chains $\gamma_{\bf T}\in Chains(FM_2)$
 such  that $Supp(\gamma_{\bf T})\subset X_{\bf T}$, where
$Supp$ means the support of a chain. In this way we obtain
a homomorphism of dg-operads $P\to Chains(FM_2)$.

4) This homomorphism is a quasi-isomorphism. This will follow
from the fact that $P$ is quasi-isomorphc to
the minimal operad  $M$, and
on the level of chain complexes (not operads) every $M_n$ is
quasi-isomorphic to $Chains(FM_2(n))$.

\begin{rmk} The theorem and its proof seem to admit a generalization to the 
case of higher dimensions.

\end{rmk}

\subsection{Little discs operad and Fulton-Macpherson operad}

We recall here the  definitions of both operads following [Ko3].

We fix the dimension $d\ge 1$.
Let us denote by $G_d$ the $(d+1)$-dimensional Lie
group acting on ${\bf R}^d$ by affine transformations $u\mapsto \lambda u +v$, 
where $\lambda>0$ is a real number and $v\in {\bf R}^d$ is a vector. 
This group acts simply transitively on the space of closed discs in
${\bf R}^d$ (in the usual Euclidean metric). The disc with center $v$ and 
with radius $\lambda$ is obtained from the standard disc 
$$D_0:=\{(x_1,\dots,x_d)\in {\bf R}^d | \,x_1^2+\dots +x_d^2\le 1\}$$
by a transformation from $G_d$ with parameters $(\lambda,v)$.
 
\begin{dfn}
The little discs operad $E_d=\{E_d(n)\}_{n\ge 0}$ 
is a topological operad defined such as follows:
 
\noindent  1) $E_d(0)=\emptyset$,
 
\noindent 2) $E_d(1)= \,\mathrm{point}\,=\{\mathrm{id}_{E_d}\}$,
 
\noindent 3) for $n\ge 2$ the space $E_d(n)$ is the space of configurations 
of $n$ disjoint discs $(D_i)_{1\le i\le n}$ inside the standard disc $D_0$.
  
The composition  $E_d(k)\times E_d(n_1)\times\dots\times E_d(n_k)\rightarrow
E_d(n_1+\dots+n_k)$ is obtained by applying elements from 
 $G_d$ associated with  discs $(D_i)_{1\le i\le k}$ in the configuration 
in $E_d(k)$ to configurations in all $E_d(n_i),\,\,i=1,\dots,k$
and putting the resulting configurations together. The action of 
the symmetric group $S_n$ on $E_d(n)$ is given by renumeration of 
indices of discs $(D_i)_{1\le i\le n}$.
\end{dfn}

The space $E_d(n)$ is homotopy equivalent to the configuration space
of $n$ pairwise distinct points in ${\bf R}^d$.

There is an obvious continuous map 
$E_d(n)\rightarrow {\mathsf{Conf}_n}(Int(D_0))$ which 
associates to a collection of disjoint discs  the collection of
their centers. This map induces a homotopy equivalence because its fibers
are contractible.

The little discs operad and homotopy equivalent little cubes operad
were introduced in topology by J. P. May 
in order to describe homotopy types of
iterated loop spaces.

The Fulton-Macpherson operad defined below is homotopy equivalent to  
the little discs operad. 
        
For $n\ge 2$ we denote by $\tilde E_d(n)$  the quotient space of
  the configuration space of $n$ points in ${\bf R}^d$
  $${\mathsf{Conf}_n}({\bf R}^d):=\{(x_1,\dots,x_n)\in ({\bf R}^d)^n| 
\,x_i\ne x_j \,\,\,{\mathrm{for\,\,\,any}}\,\,\,i\ne j\}$$
   by the action of the group $G_d$.
   The space $\tilde E_d(n)$ 
  is a smooth manifold 
of dimension $d(n-1)-1$. For $n=2$, the space $\tilde E_d(n)$ coincides with the
$(d-1)$-dimensional sphere $S^{d-1}$. There is an obvious free action of 
$S_n$ on $\tilde E_d(n)$. We define the spaces $\tilde E_d(0)$ and 
$\tilde E_d(1)$ to be 
empty. The collection of spaces $\tilde E_d(n)$ does not form an operad because 
there is no identity element, and compositions are not defined. 

Now we are ready to define the operad $FM_d=\{FM_d(n)\}_{n\ge 0}$    

The components of the operad $FM_d$ are
   
   1) $FM_d(0):=\emptyset$,
   
   2) $FM_d(1)=$point, 
   
   3) $FM_d(2)=\tilde E_d(2)=S^{d-1}$,
   
   4) for $n\ge 3$ the space $FM_d(n)$ is a manifold with corners, its
interior is  $\tilde E_d(n)$,
and all boundary strata are certain products 
of copies of $\tilde E_d(n')$ for $n^{\prime}<n$.

The spaces $FM_d(n), n\ge 2$ can be defined explicitly.
     
\begin{dfn} 
For $n\ge 2$, the manifold with corners $FM_d(n)$ is the closure
of the image of $\tilde E_d(n)$ in the compact manifold 
$\left(S^{d-1}\right)^{n(n-1)/2}\times[0,+\infty]^{n(n-1)(n-2)}$ under the map
$G_d\cdot (x_1,\dots,x_n)\mapsto
(\left({x_j-x_i\over |x_j-x_i|}\right)_{1\le i<j\le n},
{|x_i-x_j|\over|x_i-x_k|})$
where $i,j,k$ are pairwise distinct indices.

\end{dfn}
 
 One can define the natural structure of operad on the collection of 
 spaces $FM_d(n)$. We skip here the obvious definition.
      
 It is easy to check that in this way we obtain a topological operad
 (in fact an operad in the category of real compact 
 piecewise algebraic sets defined in Appendix). We call
 it the {\it Fulton-Macpherson operad} and denote by $FM_d$.
      
Set-theoretically, the operad $FM_d$ is the same as the free operad 
generated by the collection of sets $(\tilde E_d(n))_{n\ge 0}$ endowed with  
the $S_n$-actions discribed above. 

Using  piecewise
algebraic chains from  Appendix we define dg-operads
$Chains(E_d)$ and $Chains(FM_d)$ (they are operads in the
symmetric monoidal category of complexes of abelian groups). 
Since we are working over the ground field $k$ of characteristic zero,
the dg-operads will be complexes of $k$-vector spaces.

Notice that
when we write $Chains(Z)$ we mean the $cohomological$
complex with graded components $Chains_i(Z)[i], i\ge 0$ concentrated
in $negative$ degrees. 

The complex $P_n, n \ge 2$ is concentrated in degrees $[-2n-3,...,0]$.
The same is true for the complex $Chains(FM_2(n))$ (we use piecewise
algebraic chains) because $dim_{\bf R} FM_2(n)=2n-3$.

For every $n\ge 2$ the space $P_n$ has a canonical (up to signs)
basis, called the {\it standard basis}, with elements
labeled by certain combinatorial objects, which we will call
{\it meta-trees}. 

\begin{dfn} A meta-tree  ${\bf T}$ with $n$ tails  is given by the
following data:

a) an abstract tree $T\in Tree(n)$ 
together with a marking $E_i(T)\to \{finite,infinite\}$;

b) for every internal vertex $v$ of $T$ we have: $|v|\ge 2$;

c) to every $v\in V_i(T)$ we assign an
admissible labeled tree $T_v$
together with a bijection $V_{lab}(T_v)\rightarrow N_{T_v}^{-1}(v)$.

\end{dfn}

We denote by $MT(n)$ the set of isomorphism classes
of meta-trees with $n$ tails.

\begin{conj} 
There is a piecewise algebraic (see Appendix) cell decomposition
of the  spaces $FM_2(n),n\ge 2$, with the cells $\sigma_{\bf T}$ 
labeled by $MT(n)$ such that:

1) the correspondence ${\bf T}\to \sigma_{\bf T}$ is $S_n$-equivariant;

2) the operadic composition of any two cells is again a cell;

3) there is a morphism of dg-operads $Chains_{(\sigma_{\bf T})}(FM_2)\to P$
such that every cell $\sigma_{\bf T}$ is mapped (up to a sign) to
the corresponding element of the standard basis of $P$. Here
$Chains_{(\sigma_{\bf T})}(FM_2)$ denotes the chain subcomplex
of $Chains(FM_2)$ formed by $k$-linear combinations of cells $\sigma_{\bf T}$;

4) for any $n\ge 2$ the cell decomposition of $FM_2(n)$ formed by cells
$(\sigma_{\bf T}), {\bf T}\in MT(n)$ is regular, i.e. the closure
of every cell is homeomorphic to a closed ball.

\end{conj}

This conjecture implies Deligne's conjecture.

Let us now introduce a partial order on the set $MT(n), n\ge 2$.

\begin{dfn} Let ${\bf T},{\bf T}^{\prime}$ be meta-trees with $n$ tails.
We say that ${\bf T}<{\bf T}^{\prime}$ 
if there exists a sequence of meta-trees $({\bf T}_0,...,{\bf T}_m)$
such that ${\bf T}_0={\bf T}, {\bf T}_m={\bf T}^{\prime}$ 
and for any $i$ such that $0\le i\le m-1$
we have: ${\bf T}_i$ appears as a summand in the decomposition of 
$d_P{\bf T}_{i+1}$
with respect to the standard basis of $P_n$.

\end{dfn}

It follows from the condition $4)$ of the Conjecture above, that
the nerve of the partially ordered set $(MT(n),<)$ is homeomorphic
to $FM_2(n)$. The following conjecture also follows from the
Conjecture 1.

\begin{conj} The nerve of $(MT(n),<)$ is a PL-manifold with the boundary.
The above-mentioned homemorphism with $FM_2(n)$ is a homeomorphism
of PL-manifolds with boundaries.

\end{conj}

We checked this Conjecture for small $n$.

\subsection{Partial orders induced by trees}

We recall here the structure of the space $G_n$ of generators of 
the components $P_n, n\ge 1$ of the free operad $P$
(see Section 6).

There is a basis of $G_n$ elements of which are
(up to signs) parametrized by meta-trees
${\bf T}=(T,\{T_v\}_{v\in V_i(T)})$ such that
all internal edges of $T$ are finite. The degree
of the generator ${\bf T}=(T,\{T_v\}_{v\in V_i(T)})$ is equal to

$$ deg({\bf T})=\sum_{v\in V_i(T)}deg(T_v)-|E_i(T)|$$

Let $ {\bf T}$ be a generator of $P_n$.
We are going to introduce on the set $\{1,2,...,n\}$ two partial orders
$<_{h,{\bf T}}$ and $<_{v,{\bf T}}$ 
(called {\it horizontal} and {\it vertical}).
Although these orders will depend on ${\bf T}$, we will skip ${\bf T}$
from the notation if it does not lead to a confusion.

Let $i,j \in \{1,...,n\}, i\ne j$.
We have two tail vertices of $T$ labeled 
by $i$ and $j$ respectively. Then there exists a unique internal vertex $v$
of $T$ satisfying the following properties:

a) $N_T^k(i)=N_T^l(j)=v$ 
where $k,l$ are  positive integers;

b) the vertex $v$ is minimal among those satisfying a) 
(which means that $k,l$ are both minimal in a).

 Then there exist
unique labeled vertices $x,y \in T_v$ such that $N_T^{k-1}(i)=x,
N_T^{l-1}(j)=y$.
Since $T_v$ is a planar tree, we can compare $x$ and $y$ with
respect to exactly one of the following  partial orders:
``$x$ is to the left of $y$ in $T_v$" or ``$x$ is above $y$ in $T_v$".
We will call them the ``horizonal" and ``vertical" order respectively.

Let us describe the orders in $T_v$ more precisely.
 We say that $x$ is {\it above} $y$ (or $y$ is {\it below} $x$)
in $T_v$ if there exists a positive integer $a$ such that $N^a_{T_v}(y)=x$.

We say that $x$ is 
{\it to the  left} of $y$ (or $y$ is
{\it to the right} of $x$) if there exist positive integers 
$a,b$ such that $N_{T_v}^{a}(x)=N_{T_v}^{b}(y)=w$, but
$N_{T_v}^{a-1}(x)\ne N_{T_v}^{b-1}(y)$ and $N_{T_v}^{a-1}(x)$ preceeds
$N_{T_v}^{b-1}(y)$ in $N_{T_v}^{-1}(w)$ with respect to the 
order on the latter set given by the planar
structure on $T_v$.

Thus, we have defined two partial orders on  vertices of $(T_v)_{v\in V_i(T)}$.
They induce partial orders
$<_{h,{\bf T}}$ (horizontal) and $<_{v,{\bf T}}$ (vertical) on the set
$\{1,...,n\}$. We have identified the latter with the set of tails of $T$.
Namely, we say that $i<_{h,{\bf T}}j$ if (in the above notation) 
the vertex $x$ is to the left of $y$. We say that $i<_{v,{\bf T}}j$ if
$y$ is above $x$ (equivalently we say that $x$ is below $y$).

The following lemma is easy to prove.

\begin{lmm} Both horizontal and vertical orders are indeed partial
orders (i.e. they satisfy all the axioms of orders).

\end{lmm}

\begin{dfn} Let $S$ be a set, $<_1$ and $<_2$ be partial orders on $S$.
We will call them $complementary$  if any two elements
of $S$ can be compared with respect to exactly one of them.
This means that for any two elements $i,j \in S, i\ne j$ 
exactly one of the following properties holds: $i<_1j$ , $j<_1i$,
$i<_2j$, $j<_2i$.

\end{dfn}

Suppose that we have two partial orders as in the Definition.
We define on $S$ two new  pre-orders $<_{1+2}$ and $<_{1-2}$ such as
follows:

(i) $x<_{1+2}y$ if $x<_1y$ or $x<_2y$;

(ii) $x<_{1-2}y$ if $x<_1y$ or $y<_2x$.

Notice that we can reconstruct $<_1$ and $<_2$ from these
new orders. For example $x<_1y$ is equivalent to the conjunction:
 $(x<_{1+2}y)$ {\it} and $(x<_{1-2}y)$.

\begin{prp} Formulas (i) and (ii) define complete orders
on the set $S$.
\end{prp}

{\it Proof}. Straightforward.$\blacksquare$

Using this result one can easily prove the following one.

\begin{prp} Let $<_1$ and $<_2$ be  a pair of two complementary orders 
 given on a finite set $S$. Then there
 exists a unique element $s_0\in S$ such that for any $i\in S, i\ne s_0$
we have: either $s_0<_1i$ or $s_0<_2i$.

\end{prp}

$Proof$. The element $s_0$ is minimal with respect to
the complete order $<_{1+2}$.$\blacksquare$

Summarizing, we can say that 
on the set of tails of a generator of $P_n$ we have
defined two complementary partial orders (or, equivalently,
two complete orders). These orders will be used below
when we will construct the closed sets $X_{\bf T}$.

\subsection{Closed sets $X_{\bf T}$}

Let us recall that for any $n\ge 2, 1\le i,j\le n, i\ne j$ 
we have a natural projection
$p_{i,j}:FM_2(n)\to FM_2(2)$ (forgetting all  points
in $(x_1,...,x_n)$ except $x_i$
and $x_j$). As a topological space $FM_2(2)$ is identified
to the unit circle $S^1\subset {\bf R}^2$ via the map
$G_2\cdot (x_1,x_2)\mapsto {x_2-x_1\over|x_2-x_1|}$.

We denote by $S^1_{+,v}\subset FM_2$ the closed upper-half circle,
and by $S^1_{+,h}\subset S^1$ the one-element subset consisting
of the point $\{(1,0)\}$. Then $S^1_{+,v}$ corresponds to the configurations
$(x_1,x_2)\in \mathsf Conf_2({\bf R}^2)$ such that if we put 
$x_2-x_1=r e^{i \alpha}, 0\le \alpha< 2\pi$
then $\alpha\in [0,\pi]$. Similarly the
subset $S^1_{+,h}$ corresponds to the configurations $(x_1,x_2)$
such that both $x_j, j=1,2$ belong to the same horizontal line,
and $x_1$ is positioned to the left of $x_2$.

Suppose that we are given two complementary orders 
$<_h$ and $<_v$ on the set $\{1,...,n\}, n\ge 2$.
Then we define the following subset of $FM_2(n)$:

$ X_{<_h,<_v}=\{x\in FM_2(n)
|i<_hj\Rightarrow p_{i,j}(x)\in S^1_{+,h} ,
i<_vj\Rightarrow p_{i,j}(x)\in S^1_{+,v}\}.$

For a generator ${\bf T}\in P_n$ we define 
$X_{\bf T}$ as $X_{<_{h,{\bf T}},<_{v,{\bf T}}}$. 

First of all we would like to prove that $X_{\bf T}$ is contractible.
This is a special case of a more general statement.

\begin{prp} For any pair of complementary orders $<_h$ and $<_v$
given on the set $\{1,...,n\}, n\ge 2$ the subspace $X_{<_h,<_v}$
is non-empty and contractible.

\end{prp}

{\it Proof}. It can be done by
induction. For $n=2$ the result is clear, since subsets $S^1_{+,v}$ and
$S^1_{+,h}$ are contractible. Suppose that the Proposition
is true for $n-1$ points. Let us take the element $i_0$ which is
minimal in $I_n=\{1,...,n\}$ with respect to $<_{h+v}$.
We have already proved that it exists. Let us consider the induced
complementary partial orders $<_{h^{\prime}}, <_{v^{\prime}}$
on the set $I_n\setminus \{i_0\}$. They define the subset
$X_{<_{h^{\prime}},<_{v^{\prime}}}\subset FM_2(n-1)$ which is
non-empty and  contractible by induction.
Then the result is a corollary of the following observation:
the fibers of the natural projection 
$\pi:X_{<_h,<_v}\to X_{<_{h^{\prime}},<_{v^{\prime}}}$
are contractible. 
Basically it follows from the fact that the point $x_{i_0}$ is
either left or below of all the points $x_i, 1\le i\le n, i\ne i_0$.
Let us prove that the fibers of $\pi$ are non-empty.
Indeed, for any $x^{\prime}\in X_{<_{h^{\prime}},<_{v^{\prime}}}$
the operadic composition 
 $(-1,0)\circ_j x^{\prime}$ belongs to the fiber
$\pi^{-1}(x^{\prime})$. Here $(-1,0)\in S^1$ is considered as a point
in $FM_2(\{i_0,j\})$, and $j$ is an auxiliary index.

We leave to the reader the proof of contractibility of the fibers of $\pi$.
$\blacksquare$

\begin{rmk} a) One can prove that all homotopies
 can be taken in the category of
 piecewise algebraic sets.

b) It follows from the construction that 
the map ${\bf T}\mapsto X_{\bf T}$ is $S_n$-equivariant.

\end{rmk}

Now we would like to explain the Property $2$ of $X_{\bf T}$ 
stated in Section 7.1. First of all, the $S_n$-equivariance
is obvious.

Let ${\bf T}\in G$ be a generator of $P$.
We recall that  the differential in $P$ can be written schematically
(up to signs) as

$$d_P( {\bf T})=\sum_{v\in V_i(T),l}{\bf  T}_{v,l}+
\sum_{\alpha\in E_{infinite}(T)}
{\bf T}_{\alpha}^{\prime}\circ{\bf T}_{\alpha}^{\prime\prime}+
\sum_{\alpha\in E_{finite}(T),j}{\bf T}_{\alpha,j}.$$

Here meta-trees ${\bf  T}_{v,l}\in G$
 arise from the application of the differential
of $M$ to the tree $T_v$ inscribed into the vertex $v\in V_i(T)$.
This differential  was described in Section $5$. 
Index $l$ runs over all possible insertions of a new edge.

The second  summand corresponds to the tree, obtained from $T$ by
making a finite edge $\alpha$ into an infinite edge. 
The result is a composition
of  two generators of $P$ which we denote  by
${\bf T}_{\alpha}^{\prime}$ and ${\bf T}_{\alpha}^{\prime\prime}$. 

The last sum corresponds
to the operation ``contract a finite edge $\alpha$" in $T$. Then the
planar trees from $M$ inscribed into the vertices which are endpoints of 
$\alpha$
must be composed and inscribed into the
new vertex. The result is a sum of  generators ${\bf T}_{\alpha,j}$ of $P$.

\begin{prp} In the above notation we have:
$X_{{\bf  T}_{v,l}},
X_{{\bf T}_{\alpha}^{\prime}}\circ X_{{\bf T}_{\alpha}^{\prime\prime}},
 X_{{\bf T}_{\alpha,j}}$
belong to $X_{\bf T}$.

\end{prp}

{\it Proof}. Straightforward check which uses the fact that 
$S_{+,h}^1\in S_{+,v}^1$.
$\blacksquare$

\subsection{Morphism $P\to Chains(FM_2)$}

We would like to consctruct the chains $\gamma_{\bf T}\in Chains(FM_2)$ 
where ${\bf T}$ runs through the set
of generators of $P$. Let us explain the idea. 

We will construct $\gamma_{\bf T}$ by induction in the degree of $ {\bf T}$.

a) Let us assume that $deg( {\bf T})=0$. Then $T$ 
has the only internal vertex $v$ and the corresponding
planar tree $T_v$ is a binary tree. The corresponding element
of the operad $P$ is defined up to a sign. 
We make the following choice: it is just the composition
of several copies of the operation $m_2\in M_2.$

The corresponding chain $\gamma_{\bf T}$ 
should have the dimension zero. Therefore it must be a linear combination
of points with integer coefficients. We choose $\gamma_{\bf T}$ to be
just one point
in $X_{\bf T}$ with the multiplicity equal to $+1$.

b) Let us assume that $deg( {\bf T})=-1$.
 Then the formula for $d_P( {\bf T})$ contains
two summands ${\bf T}^{\prime},{\bf T}^{\prime\prime}$
corresponding to the two binary planar trees.
 It is not difficult to check that they appear with the opposite coefficients
 $\pm 1$. According to a) these summands define $0$-dimensional chains
 $\gamma_{{\bf T}^{\prime}},\gamma_{{\bf T}^{\prime\prime}}$.
 We define $\gamma_{\bf T}$ to be the only
  (up to a boundary) $1$-chain
 having
  $\gamma_{{\bf T}^{\prime}}-\gamma_{{\bf T}^{\prime\prime}}$ as the boundary.
This $1$-chain exists because of the condition imposed on multiplicities
of $0$-chains.

c) Suppose that we have a generator $ {\bf T}$,
 $deg( {\bf T})=-k, k\ge 2$.

Then we can find a chain $\gamma_{\bf T}$ of degree $-k$ such that
$$\partial(\gamma_{\bf T})=\sum_{v,l}\gamma_{{\bf T}_{v,l}}+
\sum_{\alpha} 
\gamma_{{\bf T}_{\alpha}^{\prime}\circ {\bf T}_{\alpha}^{\prime\prime}}+
\sum_{\alpha,j} \gamma_{{\bf T}_{\alpha,j}}$$
where $\partial$ is the boundary operator in the chain
complex for $FM_2$, and  rest of the notation is self-explained.

Indeed,  the RHS of this formula is known by the induction
assumption. We also know that it is a closed chain because $d_P^2=0$.
 Then one can always find a chain $\gamma_{\bf T}$ with the
given boundary. Indeed, $X_{\bf T}$ is contractible, and we 
consider chains of negative degrees (zero degree case was 
considered in a)). The space parametrizing all different choices  
of $\gamma_{\bf T}$ is contractible, 
so our choice is unique in a given homology class.

The map ${\bf T}\to \gamma_{{\bf T}}$ extends to a homomorphism
 $\Phi: P\to Chains(FM_2)$ of graded operads in such a way that
  $\Phi({\bf T})=\gamma_{{\bf T}}$.
  We have checked that
 $\Phi$ is compatible with the differentials. Therefore it is
 a homomorphism of dg-operads.

\begin{thm} The morphism $\Phi$ is a quasi-isomorphism of complexes.

\end{thm}
 
$Proof.$  We will give only a sketch of the proof. The idea is to consider a
 subcomplex $L$ of  the complex $P$ spanned by the trees
 without internal edges. This subcomplex is isomorphic 
 to the complex of the minimal dg-operad $M$.
 Indeed, $P$ can be decomposed into a direct sum
 of $L$ and a contractible complex (see Proposition $4$ in Section $4$).
 Therefore it is enough to prove that the restriction
$\Phi|_L:L\to Chains(FM_2)$ is a quasi-isomorphism of complexes.
 
 In order to do that we need to describe chains from the subcomplex
 $\Phi(L)$. Our constructions were based on certain non-canonical choices
 of chains $\gamma_{\bf T}, {\bf T}\in G$. Now, for ${\bf T}\in L$
 we will make specific choices of $\gamma_{\bf T}$.

 Let ${\bf T}\in L$ corresponds to a labeled tree  $T_v\in M$.
 We will denote it simply by $T_v$, where $T_v$ is a labeled planar
 tree inscribed in the only internal vertex $v$ of ${\bf T}$.
 
 We associate with $T_v$ an abstract tree $\widehat T_v$ such that
 $V_t(\widehat T_v)=V_i(T_v)$, and the valency of each internal vertex
 of $\widehat T_v$ is at least $2$. 
 Internal vertices of $\widehat T_v$ correspond to
 subsets $T_{\le x_0}$ described below.
 
 Let us choose an internal vertex $x_0\in V_i(T_v)$ and consider
 the set of such $x\in V_i(T_v)$ that $N_{T_v}^j(x)=x_0$ for some $j\ge 0$.
 Then internal vertices of $\widehat T_v$ correspond to such
 sets $T_{\le x_0}$ for which the cardinality $|T_{\le x_0}|\ge 2$.
 Moreover, there is a path in $\widehat T_v$
 from a tail vertex $v\in V_t(\widehat T_v)$
 to an internal vertex $T_{\le x_0}$ iff $v\in T_{\le x_0}$.
 
 It is well-known that any abstract tree $\widehat T$ with valencies
 of internal vertices at least $2$ gives rise to a stratum
 $J(\widehat T)\in FM_2(|V_t(\widehat T)|)$ (see [FM]).
 It is the operadic composition of 
 $FM_2(N^{-1}_{\widehat T_v}(u)),u\in V_i(\widehat T_v)$.
  
 We are going to
 construct a subspace $ X_{\widehat T}\subset J(\widehat T_v)$
 The space $ X_{\widehat T}$ will be constructed as the operadic
 composition of certain subspaces 
 
 $$X_u\subset FM_2(|N_{\widehat T_v}^{-1}(u)|), u\in V_i(\widehat T_v).$$
 
 Let $u=T_{\le x_0}$ be an internal vertex of $\widehat T$.
 The set of edges having  $u$ as an endpoint is in one-to-one correspondence
 with the set $\{x_0\}\sqcup N_{T_v}^{-1}(x_0).$
 The subspace $X_u$ consists of configurations of points
 $G_2\cdot (p_{x_0},(p_y)_{y\in N_T^{-1}(x_0)})$ such 
 that $p_{x_0}=(0,1)\in {\bf R}^2$,
 all points $p_y$ belong to the horizontal line 
 $\{(x,0)|x\in {\bf R}\}\subset {\bf R}^2$,
 and their order on this line is the same as their
 order in $N_{T_v}^{-1}(x_0)$.

 One can check that:
 
 a) $ X_{\widehat T_v}$ is an open cell.
 
 b) The natural projection (forgetting map)
 $FM_2(|V_i(T_v)|)\to FM_2(|V_{lab}(T_v)|)$ maps $ X_{\widehat T_v}$ onto
 an open cell $X_{T_v}\subset X_{\bf T}$. 
 
 c) the closure $\overline{X_{T_v}}$ is a manifold with
 corners (more precisely, a real piecewise algebraic
 manifold).

 d) the boundary of $\overline{X_{T_v}}$ is the union 
 of cells of the same type .

Using a)-d) one checks that the construction above gives rise
to a homomorphism of complexes $\chi:L\to Chains(FM_2)$, and
moreover, it coincides with the homomorphism $\Phi|_L$.
To be more precise, one observes that the construction of
$\Phi$ was not canonical. We made certain choices when constructed
a chain with the prescribed boundary. The point is that 
one can choose inductively the chains in such a way that 
the restriction $\Phi|_L$ of the homomorphism $\Phi:P\to Chains(FM_2)$
coincides with $\chi$.

We claim that $\chi$ induces a homotopy equivalence. This follows
from the

\begin{lmm} Let $MT_1(n)$ denotes
the subset of $MT(n)$ consisting of meta-trees with only one internal 
vertex $v$.
Then the natural embedding of the CW-complex 
$X_n=\cup_{{\bf T}\in MT_1(n)} X_{T_v}$ into $FM_2(n)$
 is a homotopy equivalence
for all $n\ge 2$.

\end{lmm}

{\it Proof}. It is not difficult to show that
$X_n$ is isomorphic to the CW-complex $\Sigma_{(n)}$
constructed in Section 5.5 via Strebel differentials.
It follows from the fact that both are regular complexes
and posets of their cells are isomorphic. Moreover, both
complexes are $K(\pi,1)$ spaces, classifying spaces for
the pure braid group of $n$-strings. Moreover, the map
$\Phi$ induces an equivalence of the fundamental groupoids.
Hence it is a homotopy equivalence. $\blacksquare$

This Lemma conludes the proof of the Theorem $2$.$\blacksquare$

\subsection{Proof of Deligne's conjecture}

The results of the previous subsection establish the Theorem 2.
In this subsection we will prove the Corollary (Deligne's conjecture).
First of all, we remark that both $P$ and $Chains(FM_2)$ are
dg-operads which are free as graded operads. Therefore the 
quasi-isomorphism morphism
of graded operads $\Phi:P\to Chains(FM_2)$ 
admits a homotopy inverse $\Psi:Chains(FM_2)\to P$ (see for ex. [M2]).

We have already proved that the Hochschild complex of an 
$A_{\infty}$-algebra is an algebra over the operad  $P$.
Indeed, $P$ is a free resolution
of $M$, and the latter dg-operad acts on the Hochschild complex.
Using the morphism $\Psi$ we make the Hochschild complex into an algebra
over the dg-operad $Chains(FM_2)$. This concludes the proof
of Deligne's conjecture.

\section{Appendix: Singular chains and differential forms}

In the Appendix we are going to describe a theory of singular
chains which is suitable for work with manifolds
with corners.
This formalism
is helpful in the proof of formality of the operad
of chains on the little disc operad (see [Ko3]).
We used it in the proof of Deligne's conjecture.

The idea of the theory of singular chains developed below
 is the following: one can construct  chains 
which produce complexes quasi-isomorphic to the standard 
complexes of singular chains,
and built out of some kind of ``piecewise algebraic spaces''.

\subsection{Spaces}

We recall that real semialgebraic sets in ${\bf R}^n$
are subsets defined by a finite number of polynomial
equations and inequalities. Constructible sets are obtained
from semialgebraic ones by boolean operations.

We define the category ${\cal P}$ of {\it compact piecewise
algebraic spaces} (compact PA-spaces for short) in the following way.

Objects of ${\cal P}$ are pairs $(X,n)$, $n=1,2,...$ such that
$X\subset {\bf R}^n$ is a compact constructible set (it is
the same as a compact real semialgebraic set).

For two objects $(X,n)$ and $(Y,m)$ the space of morphisms
$Hom((X,n),(Y,m))$ is formed by continuous maps $f:X\to Y$
such that  $graph(f)\subset {\bf R}^n\times {\bf R}^m$ is
constructible.

In the future we will
 skip the index $n$ in the notation
$(X,n)$ if it will not lead to a confusion.

Obviously we have a functor from ${\cal P}$ to the category of
compact Hausdorff topological spaces. An isomorphism in
${\cal P}$ is a morphism which is a homeomorphism of topological
spaces.

Let $X\in {\cal P}$. Then one can define a sheaf ${\cal O}_X$
of piecewise algebraic functions on $X$. To do this we note
first that one can speak about constructible subsets of $X$.
By definition they are constructible sets in the bigger space
${\bf R}^n$. Then for any open $U\subset X$ we define ${\cal O}_X$
to be an ${\bf R}$-algebra of continuous functions $f$ on $U$ 
such that for any compact constructible $V\subset U$ we have:
$graph (f|_V)$ is constructible. It is easy to see that we get 
a sheaf ${\cal O}_X$ of algebras.

\begin{lmm} Morphisms from $X$ to $Y$ in the category ${\cal P}$
are in one-to-one correspondence with homomorphisms
  of algebras
${\cal O}(Y)\to {\cal O}(X)$.

\end{lmm}

{\it Proof}. Exercise. $\blacksquare$

Similarly, we can consider non-compact case.

\begin{dfn} Piecewise algebraic space (PA-space for short) is
a locally compact Hausdorff topological space $X$, equipped with
the sheaf ${\cal O}_X$ of ${\bf R}$-algebras which is locally
isomorphic to ${\cal O}_{X^{\prime}}$ for some compact 
$X^{\prime}\subset {\cal P}$.

\end{dfn}

Clearly PA-spaces form a category. Compact PA-spaces are exactly
objects of ${\cal P}$.

We define a $d$-dimensional PA-manifold with boundary as a PA-space $X$
which is modeled locally by the closed half-space ${\bf R}^d_+$.

\subsection{Singular chains}

We would like to define an appropriate version of singular chains
for  PA-spaces. We will give two equivalent descriptions.

{\it First description}.

Let $X$ be a PA-space. We define the group of 
 $n$-chains $Chains_n(X):=C_n(X,{\bf Z})$ as an abelian group
generated by equivalence classes of triples $(M,or,f)$ such that:

i) $M$ is a compact PA-manifold of dimension $n$;

ii) $f:M\to X$ is a morphism in ${\cal P}$;

iii) $or$ is an orientation of $M$.

We need to define an equivalence relation. It is the same as to say when
a finite linear combination $\sum_in_i(M_i,or_i,f_i)$ is equal to zero
in $C_n(X,{\bf Z})$.

Notice that $Y=\cup_if_i(M_i)$ carries the structure of a compact
 PA-space, of the dimension $\le n$. Then there exists a constructible
subset $Y_0\subset Y$ such that $dim(Y\setminus  Y_0)\le n-1$, and  for
any point $y\in Y_0$, any $i$ and $x_{i,\alpha}\in f_i^{-1}(y)$ there
exists a neighborhood $U_{i,\alpha}$ of $x_{i,\alpha}$ such that 
$Y_0\cap f_i(U_{i,\alpha})$
 is a PA-manifold and the
morphism $f_i|_{U_{i,\alpha}}$ is a
 homeomorphism of $U_{i,\alpha}$ onto its image.

We choose an orientation $or_y$ of $Y_0$ near the point $y$. Then the above-
mentioned linear combination is declared to be  zero iff for every 
point $y\in Y_0$ we have $\sum_{i,\alpha}n_isgn(i,\alpha,y)=0$ where
$sgn(i,\alpha,y)$ is defined as 
$f_{i*}|_{U_{i,\alpha}}(or_i)=sgn(i,\alpha,y)or_y$,
and sum is taken over all points $x_{i,\alpha}$ such that
$f_i(x_{i,\alpha})=y$. Notice that $sgn(i,\alpha,y)$ always takes
values $\pm 1$.

{\it Second description}.

We will define $C_n(X)=C_n(X,{\bf Z})$ as a quotient by certain
equivalence relation of the set
of quadruples
$(Y,Y_0,or,mult)$ such that:

a) $Y\subset X$ is a compact 
PA-subspace, $dim Y\le n$, which contains an open dense
constructible subset $Y_0$ without boundary 
which is a PA-manifold  of dimension $n$;

b) $or$ is an orientation of $Y_0$;

c) $mult:Y_0\to \{1,2,...,\}$ is a locally constant map
(multiplicity).

The equivalence relation is defined such as follows.

We say that  $(Y,Y_0,or,mult)$ is equivalent to $(Y^{\prime},
Y_0^{\prime},or^{\prime},mult^{\prime})$ iff

(i) $Y=Y^{\prime}$ and there exists an open contsructible
$Y_0^{\prime \prime}\subset Y_0^{\prime}\cap Y_0^{\prime \prime}$
such that $Y$ is equal to the closure of $Y_0^{\prime \prime}$;

(ii) restrictions of orientations and
multiplicity functions of $Y_0$ and $Y_0^{\prime}$
to $Y_0^{\prime\prime}$ coincide.

Using the pairing between chains and differential
forms (see 8.3), one can show that the first and the second 
descriptions give canonically isomorphic sets of chains.

One can also check the
 following properties of singular chains:

1) there is a naturally defined differential $\partial:
C_n(X,{\bf Z})\to C_{n-1}(X,{\bf Z}),\partial^2=0$ (in the Description 1
it is well-defined by the formula $\partial(M,or,f)=
(\partial M,\partial (or), f|_{\partial M})$).

2) the correspondence $X\to C_{\cdot}(X,{\bf Z})$
is a functor from the category of PA-spaces to the category
of abelian groups;

3) $C_i(X,{\bf Z})$ vanishes for $i<0$ and $i>dim X$;

4) if $X$ is a compact oriented PA-manifold, then
there exists a canonically defined chain $[X]\in  C_n(X,{\bf Z})$.
In the Description 1 it is defined by the formula
$[X]=(X,or,1)$;

5) for any finite
collection of PA-spaces $(X_i)_{i\in I}$
 there is a natural homomorphism of complexes of abelian groups
$\otimes_{i\in I}C_{\cdot}(X_i,{\bf Z})\to 
C_{\cdot}(\prod_{i \in I} X_i,{\bf Z})$;

6) there is a naturally defined soft sheaf
of complexes $\underline C_{\cdot}^{closed}$
on $X$ such that for a compact $X$ the abelian group
$\Gamma(X,\underline C_{\cdot}^{closed})$ coincides
with $C_{\cdot}(X,{\bf Z})$.

7) the homology of $C_{\cdot}(X,{\bf Z})$ is naturally isomorphic to the usual
singular homology $H_{\cdot}(X,{\bf Z})$.

For non-compact PA-spaces we define locally finite chains
as global sections of the sheaf $\underline C_{\cdot}^{closed}$.
Then the homology 
of $C_{\cdot}^{closed}(X,{\bf Z})$ is isomorphic to the usual
singular homology with locally compact support
$H_{\cdot}^{closed}(X,{\bf Z})$.

\begin{rmk}

a) Mayer-Vietoris sequence trivializes for PA chains:
if $X=Y_1\cup Y_2$ is a union of locally closed constructible
subsets, then 
 the corresponding short sequence of abelian groups
$0\to C_n(Y_1\cap Y_2,{\bf Z}) \to C_n(Y_1,{\bf Z})\oplus C_n(Y_2,{\bf Z})
\to C_n(X,{\bf Z})\to 0$
 is exact.

 b) The Property 5 is very convenient in order to formulate
Deligne's conjecture. It allows to avoid Eilenberg-Zilber theorem.

c) Property 4 seems to be useful for the higher-dimensional generalization
of Deligne's conjecture (see [Ko3]). Namely, if $A$ is a $d$-algebra (i.e. 
an algebra over the dg-operad $Chains(FM_d)$), then $A[d-1]$
carries a natural structure of an $L_{\infty}$-algebra.
More precisely, for any $n\ge 2$
the fundamental cycle $[FM_d(n)]$ gives the $n$-th higher bracket
on $A[d-1]$.

\end{rmk}

\subsection{Differential forms: first approximation}

Let $X$ be a PA-space.  We would like to define
an appropriate notion of differential form on $X$. In this
subsection we construct a first approximation to the future
algebra of differential forms. 

\begin{dfn} Sheaf  $\Omega^k_{X,min}$ is locally defined
as a vector subsubspace in $Hom(C_k(X,{\bf Z}),{\bf R})$ generated
by functionals $l=(f_0,f_1,...,f_k), f_i\in {\cal O}_X$ such that

$$l(M,or,\phi)=\int_{M_0}\phi^{\ast}(f_0df_1\wedge df_2...\wedge df_k).$$

\end{dfn}

Here $M_0=M$ is a dense open constructible subset.
We can assume that $M$ is a subspace of some ${\bf R}^N$,
and all functions $f_i$ are smooth on a smooth dense open submanifold
$M_0\subset M$. 

\begin{prp} Functional $l$ from the definition above is well-defined
(i.e. the integral absolutely converges).

\end{prp}

{\it Proof}. We can use a sequence of functions $(f_0,f_1,...,f_k)$
in order to map $M$ to ${\bf R}^{k+1}$. Then we obtain a
singular chain $\gamma \in C_k({\bf R}^{k+1})$.
The support of $\gamma$ is a compact constructible subset of
dimension $\le k$. Moreover, the volume
of $\gamma$ with respect to the induced metric is finite.
It follows that $\int_{\gamma}\omega$  absolutely
converges for an arbitrary smooth differential
$k$-form $\omega$. In particular converges the integral
in question.$\blacksquare$

It is easy to check that the differential $d:\Omega^k_{X,min}\to
\Omega^{k+1}_{X,min}$ is well-defined (as the adjoint to the 
boundary operator on the chains). Moreover the following 
Stokes formula holds.

\begin{lmm} In the previous notation we have:
$$(1,f_0,...,f_k)(M,or, \phi)=
(f_0,...,f_k)(\partial M,or|_{\partial M},\phi|_{\partial M})$$

\end{lmm}

Hence we obtain a complex of sheaves $\Omega^k_{X,min}$ for every
$k\ge 0$. Our (first approximation to) $k$-forms are sections 
of these sheaves. One can take pull-backs of forms.

Moreover, we can introduce wedge product of forms. Indeed the space
of $k$-forms $\Omega^k_{min}(X)$ is naturally embedded
into the direct limit of spaces of smooth $k$-forms
 $\Omega^k_{C^{\infty}}(X_{\alpha})$. Here $X_{\alpha}\subset X$
 are open dense locally constructible subspaces,
 $X_{\alpha}=\sqcup_i X_{i,\alpha}$, where $X_{i,\alpha}$
 is a PA-manifold. Moreover, it is assumed that we have fixed
 an embedding $X_{i,\alpha}\to {\bf R}^{N_{i,\alpha}}$ identifying
 $X_{i,\alpha}$ with a $C^{\infty}$-submanifolds
 in ${\bf R}^{N_{i,\alpha}}$. Using the embeddings
 we define $\Omega^k_{C^{\infty}}(X_{\alpha})$.
 Therefore  the wedge product given by the formula
$(f_0,...,f_k)\wedge (g_0,...,g_n)=
 (f_0g_0,f_1,...,f_k,g_1,...,g_n)$ 
 is well-defined.

The above-defined forms 
are sections of soft sheaves because
one can use piecewise algebraic functions to produce
a partition of unity, so the standard proofs work.
One can show that Mayer-Vietoris sequence degenerates as in the case
with chains above: if $X=Y_1\cup Y_2$ is a union of locally-closed
piece-wise algebraic subspaces then one has a short exact sequence
$$ 0\to \Omega^k_{min}(X)\to \Omega^k_{min}(Y_1)\oplus \Omega^k_{min}(Y_2)\to
\Omega^k_{min}(Y_1\cap Y_2)\to 0
$$

In the next subsection we will extend algebras $\Omega^{\ast}_{min}(X)$
adding push-forwards of such forms. The following example illustrates
one of the reasons for that.

\begin{exa} Let $Y=[0,1]\times [0,1]$, $X=[0,1]$ and $f:Y\to X$
be the map $(x,y)\mapsto t=xy$. Then, outside of $t=1$,
 we have a bundle in the category
${\cal P}$ (i.e. base, fibers and total space are objects,
 projection is a morphism
in ${\cal P}$). Take the form $\omega=xdy$. Then $f_{\ast}\omega=t(logt)$
is a continuous function, but does not belong to ${\cal O}(X)$.

\end{exa}

\subsection{Full algebra of forms}

We are going to define a dg-algebra $\Omega_{PA}^{\cdot}(X)$ which
contains $\Omega^{\cdot}_{min}(X)$, satisfies Poincare lemma, etc.
Elements of $\Omega_{PA}^{\cdot}(X)$ will be called PA-forms,
although their coefficients are not PA-functions (see the previous
Example).
It is interesting that only $0$-forms will
 be forms with continuous coefficients. Higher
order forms can have coefficients in $L_1^{loc}(X)$.
But they are still closed
under the wedge product.

We start with some preliminaries. Let $f:Y\to X$ be a 
proper 
PA-map of PA-spaces. For simplicity we will give
the definition below in the case when $X$ is compact.

\begin{dfn} Continuous family of PA $k$-cycles (PA-family
of cycles for short) is defined as an element of the abelian 
group $C_k(Y\to X)$ described below. 

The group $C_k(Y\to X)$ consists of maps $\Phi:X\to C_k(X)$ such that:

a) For any $x\in X$ the set $Supp\,\Phi(x)$ belongs to $f^{-1}(x)$.

b) The set $Z=\cup_{x\in X}Supp\, \Phi(x)$ is constructible with the compact
closure  $\overline {Z} \subset Y$.

c) There exists a dense constructible $Z_0\subset \overline {Z}$ such that
for any $x\in X$ the intersection $Z_0\cap f^{-1}(x)$ is a
PA-manifold without boundary,
 dense in $\overline {Z}\cap f^{-1}(x)$. Moreover,
the chain $\Phi(x)$ is obtained from some orientation of $Z_0\cap f^{-1}(x)$
and locally constant multiplicity map $Z_0\cap f^{-1}(x)\to {\bf Z}$ 
as in the second description  of PA-chains.

d) For any $x\in X$ and $z\in f^{-1}(x)\cap Z_0$, and any sequence $(x_i),
x_i\to x$ 
 the multiplicity at the
point $z$ is the ``natural limit'' of the multiplicities of 
$f^{-1}(x_i)\cap Z_0$.

\end{dfn}

In the example below $X$ not compact. It easy to make necessary
modifications of our definition, so it will be valid in the non-compact case.
In that case we get a sheaf on $X$.

\begin{exa} Let $Y={\bf R}^2, X={\bf R}$ and $f:Y\to X$ is the projection
$(x_1,x_2)\mapsto x_1$. Let us define $\Phi(x_1)$ to be equal 
$\{x_1\}+\{-x_1\}$ for $x_1>0$ and equal to $2\{0\}$ otherwise.
Then this map gives an element of $C_0(Y\to X)$.

\end{exa}

We will briefly describe some properties of PA-families of cycles below.

1) PA-families $C_{\cdot}(Y\to X)=\{C_k(Y\to X)\}_{k\ge 0}$ 
form a complex of sheaves (abelian groups in the case when $X$
is compact).

2) There are natural operations $C_k(Y\to X)\otimes C_l(X)\to
C_{k+l}(Y)$. Informally we will denote them by $(\gamma,\alpha)\mapsto 
\gamma \times \alpha$.

3) For any Cartesian square 

$$\begin{array}{lr}
Y^{\prime}& \to Y\\
\downarrow & \downarrow\\
X^{\prime}& \to X\\
\end{array}
$$
one has the natural morphisms
$$C_k(Y\to X)\to C_k(Y^{\prime}\to X^{\prime})$$
and
$$C_k(Y\to X)\otimes C_l(X^{\prime}\to X)\to C_{k+l}(Y^{\prime}\to X).$$

Now we are ready to define piecewise PA-forms.

\begin{dfn}
PA-forms of degree $k$ are sections of soft sheaves which are locally given
by functionals $l:C_k(X)\to {\bf R}$ such that
$$l(\alpha)=\int_{\gamma\times \alpha}\omega$$
where $\gamma \in C_l(Y\to X),l\ge 0$ for some proper $Y\to X$,
and $\omega \in \Omega^{k+l}(Y)_{min}$.

\end{dfn}
 
 We will denote the space of  
 PA-forms on $X$ by $\Omega_{PA}^{\cdot}(X)$.

Similarly to the case of $\Omega_{min}^{\cdot}(X)$ one can see that
at the ``generic point'' the space $\Omega_{PA}^{\cdot}(X)$ can be
naturally embedded into the space of differential forms.
In fact we obtain a soft sheaf of dg-algebras. It is closed 
under pushforwards $f_{\ast}$ where $f:Y\to X$ is a locally trivial
bundle  in the category ${\cal P}$ with fibers which are
compact oriented PA-manifolds.

The latter fact can be generalized further. One can consider
a family $f:Y\to X$ where $Y$ and $X$ are oriented
compact PA-manifolds, and all fibers of $f$ have the same
dimension $dim(Y)-dim(X)$.
 Then one get a continuous family
of cycles (``fundamental cycle'' of $f^{-1}(x), x\in X$)
over $X$. This gives a pushforward of PA-forms for certain maps
(``flat morphisms") which are not necessarily fibrations.

Notice that the Poincare lemma holds for piecewise forms. One
can imitate the usual proof based on the integration over rays
$(x,t)\mapsto xt$, where $x\in U\subset {\bf R}^n$, $U$ is a convex
domain containing $0\in {\bf R}^n$, and $t\in [0,1]$.

\vspace{15mm}

{\bf REFERENCES}

\vspace{5mm}

[BD] A. Beilinson, V. Drinfeld, Chiral algebras, preprint (1995).

\vspace{2mm}

[B] R. Borcherds, Vertex algebras, q-alg/9706008.

\vspace{2mm}

[BV] J. Boardman, R. Vogt, Homotopy invariant algebraic structures on
topological spaces, Lect. Notes Math. 347 (1973).

\vspace{2mm}

[De] P. Deligne, letter to Drinfeld (1990)

\vspace{2mm}

[Dr] V. Drinfeld, On quasi-triangular quasi-Hopf algebras and a group
closely related with $Gal(\overline{\bf Q}/{\bf Q})$, Leningrad Math. J
vol. 2 (1991), 829-860.

\vspace{2mm}

[FM] W. Fulton, R. MacPherson, A compactification of configuration spaces,
Annals Math., 139(1994), 183-225.

\vspace{2mm}

[GeJ] E. Getzler, J. Jones, Operads, homotopy algebra, and iterated integrals
for double loop spaces. Preprint 1994, hep-th/9403055.

\vspace{2mm}

[GeKa] E. Getzler, M. Kapranov, Cyclic operads and cyclic homology.
In: Geometry, Topology  and Physics for Raul Bott. International Press,
1995, p. 167-201.

\vspace{2mm}

[GiKa] V. Ginzburg, M. Kapranov, Koszul duality for operads,
Duke Math. J. 76 (1994), 203-272.

\vspace{2mm}

[KM] M. Kontsevich, Yu. Manin, Gromov-Witten classes, quantum cohomology,
and enumerative geometry, Comm. Math. Phys., 164:3 (1994), 525-562.

\vspace{2mm}

[Ko1] M. Kontsevich, Deformation quantization of Poisson manifolds,I.
Preprint q-alg/9709040.

\vspace{2mm}

[Ko2] M. Kontsevich, Rozansky-Witten invariants via formal geometry.
Preprint q-alg/9704009.

\vspace{2mm}

[Ko3] M.Kontsevich, Operads and motives in deformation quantization.
Letters in Math.Phys., 48:1 (1999). Preprint math.QA/9904055

\vspace{2mm}

[M] M. Markl, Cotangent cohomology of a category and deformations.
J. Pure Appl. Alg., 113:2 (1996), 195-218.

\vspace{2mm}

[M2] M. Markl, Homotopy Algebras are Homotopy Algebras. Preprint
  math.AT/9907138
   
\vspace{2mm}

[Ma] J. May, Infinite loop space theory, Bull. AMS, 83:4 (1977), 456-494.

\vspace{2mm}

[MS] J. McClure, J. Smith, A solution of Deligne's conjecture, math.QA/9910126.

\vspace{2mm}

[So] Y. Soibelman, Meromorphic tensor categories, q-alg/9709030. To
appear in Contemporary Math.

\vspace{2mm}

[St] K. Strebel, Quadratic differentials. Berlin: Springer-Verlag, 1984.

\vspace{2mm}

[T] D. Tamarkin, Formality of chain operad of small squares, math.QA/9809164.

[V] A. Voronov, Homotopy Gerstenhaber algebras, math.QA/9908040

\vspace{20mm}

Authors addresses: I.H.E.S., 35 route de Chartres,
 Bures-sur-Yvette 91440, France (M.K. and Y.S.), and Department of Mathematics,
 Kansas State University, Manhattan, KS 66506, USA (Y.S.).
 
\vspace{5mm} 
 emails: maxim@ihes.fr, soibel@ihes.fr

\end{document}